\documentclass[12pt]{article}
\usepackage{epsfig}
\usepackage{amssymb}
\usepackage{amsfonts}
\usepackage{graphicx}
\usepackage{epstopdf}

\newtheorem{theorem}{Theorem}[section]
\newtheorem{lemma}[theorem]{Lemma}
\newtheorem{corollary}[theorem]{Corollary}
\newtheorem{proposition}[theorem]{Proposition}

\newtheorem{remark}[theorem]{Remark}

\newtheorem{definition}[theorem]{Definition}
\newtheorem{example}[theorem]{Example}

\newtheorem{exercise}[theorem]{Exercise}

\newenvironment{proof}{{\it Proof:\/}}{$\Box$\vskip 0.08in}
\newcommand{\mod}{{\mbox{ mod }}}

\newcommand\Z{{\mathbb Z}}
\newcommand{\lk}{{\mbox{ lk }}}

\begin{document}
\centerline{\bf The Trieste look at Knot Theory;}
\centerline{J\'ozef H. Przytycki (Washington)}
\ \\
{\bf Abstract}
{\footnotesize
This paper is base on talks which I gave in May, 2010 at Workshop in Trieste (ICTP).
In the first part we present an introduction to knots and knot theory 
from an historical perspective, starting from Summerian knots and ending 
on Fox 3-coloring. We show also a relation between 3-colorings and the Jones polynomial. 
In the second part we develop the general theory of Fox colorings and 
show how to associate a symplectic structure to a tangle boundary so that tangles 
becomes Lagrangians (a proof of this result has not been published before). 
We also discuss rational moves on links and their relation to Fox colorings. 
}

\section{ Classical Roots of Knot Theory}
Knots have fascinated people from the dawn of human history. One of the oldest 
examples of knots in art or religion 
is the cylinder seal impression (c. 2600-2500 B.C.)
 from Ur, Mesopotamia (see Figure 1.1). It is  described
in the book {\it Innana} by Diane Wolkstein and
Samuel Noah Kramer \cite{Wo-Kr} (page 7),
illustrating the text:

``Then a serpent who could not be charmed
made its nest in the roots of the tree."
\ \\
\ \\
\centerline{\ \ \ \ \ \ \ \ \ \ \ \ \
\psfig{figure=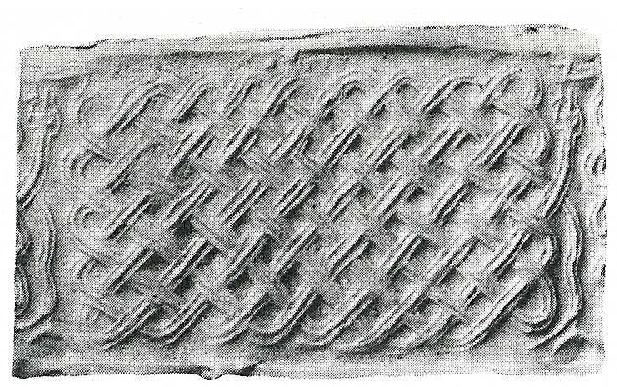,height=6.9cm}}
\centerline{Fig. 1.1; \ {\it Snake with Interlacing Coil}}
{\footnotesize {\it
Cylinder seal. Ur, Mesopotamia. The Royal Cemetery, Early Dynastic
period, c. 2600-2500 B.C. Lapis lazuli. Iraq Museum. Photograph
courtesy of the British Museum, UI 9080,
\cite{Wo-Kr}.}}

In today's audience we have students from all over the world, including Iraq, in which Ur lies today. 
I encourage you to send me an example of an ancient knot 
from your culture or country.

In the 19'th century knot theory was an experimental science.
Topology (or geometria situs) had not developed enough to
offer tools that allow precise definitions and proofs.
Johann Benedict Listing (1808-1882), a student of Gauss and pioneer of knot theory, writes in \cite{Lis}:
\emph{In order to reach the level of exact science, topology will
have to translate facts of spatial contemplation into easier notion
which, using corresponding symbols analogous to mathematical ones, we
will be able to do corresponding operations following some simple rules.} 
A combinatorial definition of the Gaussian linking number (initially defined by Gauss in 1833 as an integral \cite{Gaus}) 
was the first step in realizing Listing's program, \cite{Brunn-1892}. 

Much of early knot theory was motivated by physics and chemistry. In the
1860s, it was believed that a substance called the ether pervaded all of space. In
an attempt to explain the different types of matter, William
Thomson later known as  Lord Kelvin (1824-1907)  hypothesized that atoms were merely knots in the fabric of this
ether. Different knots would then correspond to different elements.
 
Thus, in the second half of the 19'th century knot theory was developed
primarily by physicists (Thomson, James Clerk Maxwell (1831-1879),
Peter Guthrie Tait(1831-1901)), and one can argue that a high
level of precision was not always 
appreciated\footnote{For an outline of the global history of knot theory see  
\cite{P-4} or the second chapter
of my book on knot theory \cite{P-Book}.}.

\subsection{Felix Klein observation on Knotting in dimension four}\ \\
 Tait write in his paper of  1877 \cite{Tait-2}: {\it Klein himself
made the very singular discovery that in space of four
dimensions there cannot be knots. Klein observation was noticed
in non-mathematical circles and it became part of popular culture.
For example, the American magician and medium
Henry Slade was performing ``magic tricks" claiming that he
solves knots in fourth dimension. He was taken seriously by a German
astrophysicist J.K.F.Zoellner who had with him a number of seances
in 1877 and 1878.}.
Tait is referring to ''Mathematische Annalen, ix. 478" and later authors often cite 
that paper (\cite{Klein}, 1875), and even give the page 476 (\cite{As}). 
Most likely it is a misunderstanding as Klein  discusses
 there intrinsic and ambient topological properties (of curves and surfaces) but never in 
context of knots in dimension four. More likely explanation is that Klein described to 
Tait the observation in a private correspondence. For many of you German is a native language.
I would challenge you to go through Klein papers and correspondence to find a root of 
Klein's ``singular discovery". 
\\

\subsection{Precision comes to Knot Theory}
Throughout the 19'th century knots were understood as closed curves in space
up to a natural deformation, described as a movement in
space without cutting and pasting. This understanding allowed scientists
(Tait, Thomas Penyngton Kirkman, Charles Newton Little, Mary Gertrude Haseman)
to build tables of knots, but did not lead to precise methods
  to distinguish knots that cannot be be practically deformed into each other.
In a letter to O.~Veblen, written in 1919, young J.~Alexander expressed
his disappointment\footnote{We should remember that this was written by
a young revolutionary mathematician forgetting that he was ``standing on
the shoulders of giants.''. In fact, the knot invariant Alexander outlined in the letter
is closely related to the Kirchhoff matrix, and the numerical invariant he also obtained 
is equivalent to complexity of the signed graph corresponding to the link via Tait translation;
see Subsection 1.3.}:
\emph{``When looking over Tait \emph{On Knots} among other things,
He really doesn't get very far. He merely writes down all the plane projections
of knots with a limited number of crossings, tries out a few transformations
that he happens to think of and assumes without proof that if he is unable
to reduce one knot to another with a reasonable number of tries, the two
are distinct. His invariant, the generalization of the Gaussian invariant \dots
for links is an invariant merely of the particular projection of the knot
that you are dealing with, - the very thing I kept running up against
in trying to get an integral that would apply."} 

In 1907, in the famous {\it Mathematical Encyclopedia}, Max Dehn and
Poul Heegaard outlined a systematic approach to topology.
In particular, they precisely formulated the subject of
 knot theory \cite{D-H}. To bypass the notion
of deformation of a curve in a space (which was not well defined at the time)
they introduced lattice knots and a precise definition of
 (lattice) equivalence.
Later on, Reidemeister and Alexander considered more general
polygonal knots in a space, with equivalent knots related by
a sequence of $\Delta$-moves; they also explained $\Delta$-moves via elementary moves
on link diagrams -- Reidemeister moves (see Figure 1.6).
 The approach of Dehn and Heegaard
was long ignored, however recently there has been interest in the study of lattice knots\footnote{I am aware 
of two exceptions: in 1954, a popular article of Alan Turing (1912-1954) considers elementary 
moves on knots that lie on the unit lattice in $R^3$. He concludes: ``A similar decision problem 
which might well be unsolvable is the one concerning knots which has already been mentioned." 
\cite{Turing},\cite{Gor-2}. In 1962, the biophysicist Max Delbr\"uck (1906 -- 1981) winner of 
the Nobel Prize in Physiology or Medicine in 1969,
proposed that long molecules discovered in living organisms can be knotted, and asked about the shortest
length of such a knot \cite{Del}. In his model, lattice knots are restricted to those having straight
segments of length 1. Delbr\"uck found a realization of the trefoil knot of length 36.
Delbr\"uck's problem was popularized by Martin Gardner (1914-2010) in the November 1970 issue of Scientific American,
where Gardner had a popular ``Mathematical Games" column.
Gardner comments that it is still unknown whether 36 is the minimal number of segments for
Delbr\"uck's  (molecule) lattice nontrivial knot and he comments that if a segment can be of any length,
then 24 is possible. We know now that
this is the smallest number \cite{Dia}.} (e.\,g.\ \cite{B-L}).

\subsection{Early invariants of links}
The fundamental problem in knot theory was, until recently,\footnote{There are
now algorithms that allow recognition of any knot, but they are very slow \cite{Mat}. 
Modern knot theory, on the other hand,  looks for structures on a space of knots or for 
mathematical or physical meanings of knot invariants.}
to distinguish
non-equivalent knots. Even in the case of the unknot and the trefoil knot,
this was not achieved until the fundamental work of
Jules Henri Poincar\'e (1854-1912) was applied. In his seminal paper ``Analysis Situs''
 (\cite{Po-1} 1895) he laid the foundations for algebraic topology.
According to W.~Magnus \cite{Mag}:
\emph{Today,  it appears to be a hopeless task to assign priorities for the
definition and the use of fundamental groups in the study of knots,
particularly since Dehn had announced \cite{De-0} one of the important
results of his 1910 paper (the construction of Poincar\'e spaces with the
help of knots) already in 1907}. Wilhelm Wirtinger (1865-1945), in his lecture
delivered at a meeting of the German Mathematical Society in 1905 outlined
a method of finding a knot group presentation (now called the
Wirtinger presentation of a knot group, 
but examples using his method
were given only after the work of Dehn.

\subsection{Tait's relation between knots and graphs.}\label{Subsection 1.3}

   Tait was the first to notice the relation between knots and
planar graphs. He colored the regions of a knot diagram alternately white
and black so that the infinite region is black. He then constructed a graph by placing a vertex
inside each white region, and connecting vertices by edges going through
the crossing points of the diagram (see Figure 1.2).


\ \\
\centerline{\psfig{figure=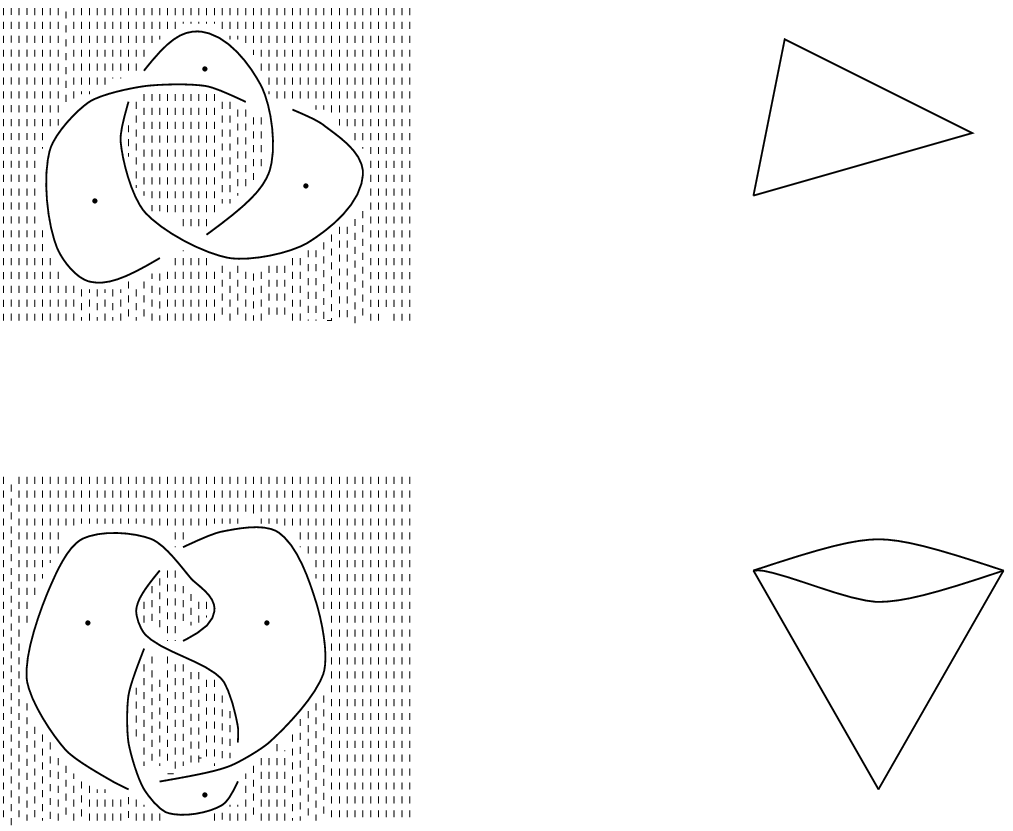,height=8.5cm}}\ \\

\centerline{Figure 1.2; Tait's construction of graphs from link diagrams as described in \cite{D-H}}\ \\
It is useful to mention the Tait construction in the opposite direction, going from
a signed planar graph $G$ to a link diagram $D(G)$. We replace every signed edge of a graph by a crossing
according to the convention of Figure 1.3, and connect endpoints along edges as in
Figures 1.4 and 1.5.

\centerline{\psfig{figure=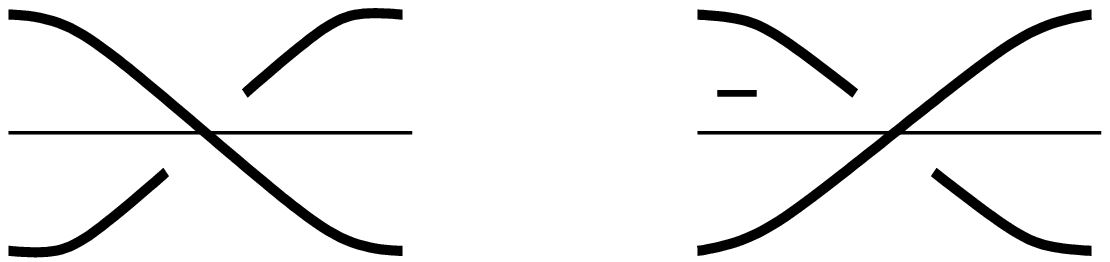,height=2.5cm}}\ \\
\centerline{Figure 1.3; convention for crossings of signed edges (edges without markers \\
 are assumed to be positive)}\ \\

\centerline{\psfig{figure=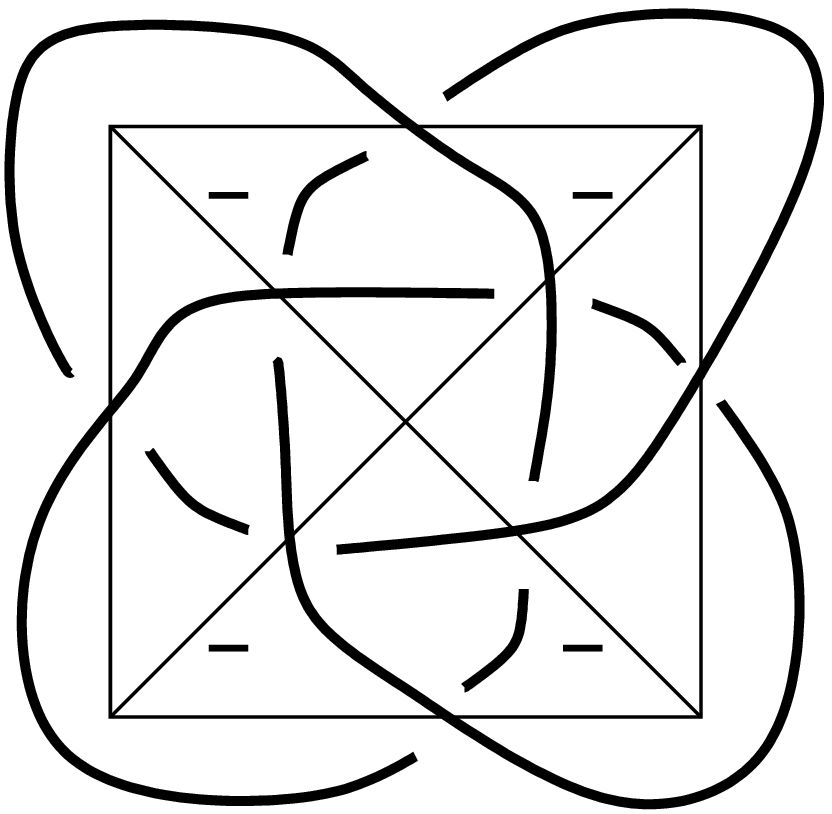,height=3.7cm}}\ \\

\centerline{Figure 1.4; The knot $8_{19}$ and its Tait graph ($8_{19}$ is the first non-alternating knot in tables)}\ \\


We should mention here one important observation already known to Tait (and in explicit form to Listing):
\begin{proposition}\label{Proposition 1.6}
The diagram $D(G)$ of a connected graph $G$
is alternating if and only if $G$ is positive (i.\,e.\ all edges of $G$ are positive) or
$G$ is negative.
\end{proposition}
A proof is illustrated in Figure 1.5. \ \\

\centerline{\psfig{figure=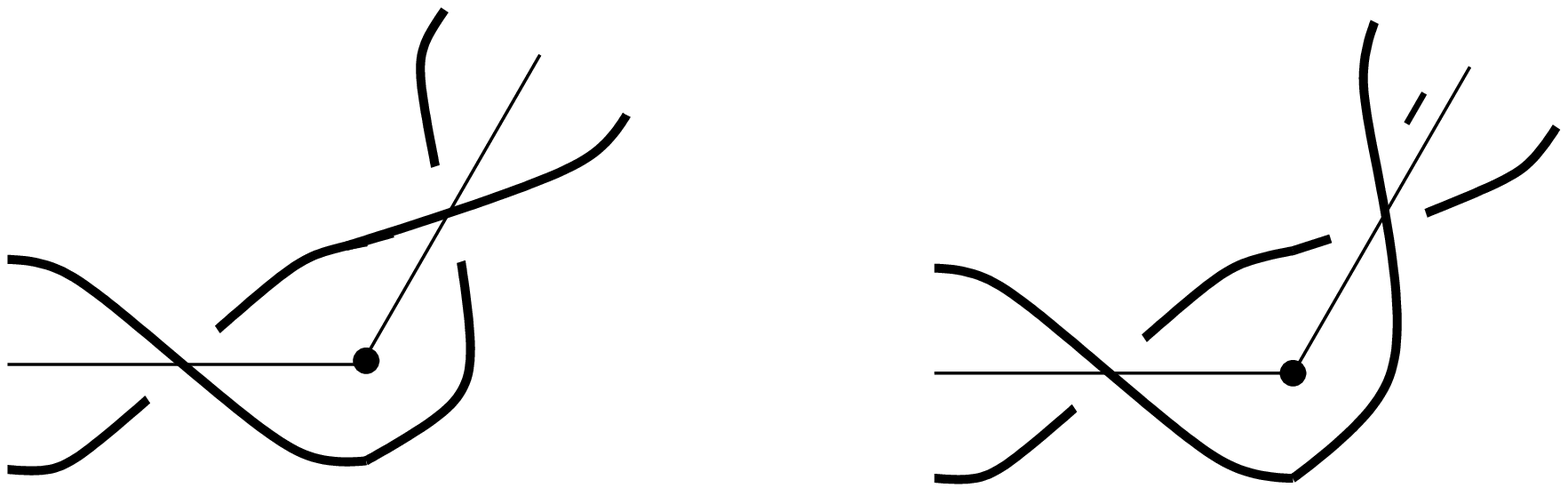,height=3.1cm}}\ \\
\centerline{Figure 1.5; Alternating and  non-alternating parts of a diagram}
\begin{exercise} Draw all connected plane graphs of up to 7 edges without loops and isthmuses (edges 
whose removal disconnects a graph). Identify related Tait diagrams with knots and links in tables 
of knots \cite{Rol}.
\end{exercise}

Maxwell was the first person to consider the question of two projections representing
equivalent knots. He considered some elementary moves (reminiscent of the future Reidemeister moves),
but never published his findings.

The formal interpretation of equivalence of knots in terms of diagrams
was described by Reidemeister \cite{Rei-1}, 1927, and Alexander and Briggs
\cite{A-B}, 1927.
\begin{theorem} [Reidemeister theorem] \label{Theorem 1.3}

Two link diagrams are equivalent\footnote{In modern knot theory, especially after the
work of R.~Fox, we use
usually the equivalent notion of ambient isotopy in $R^3$ or $S^3$. }
 if and only if they are connected by a
finite sequence of Reidemeister moves $R_i^{\pm 1}, i=1,2,3$ (see Fig.~1.6)
and isotopy of the  diagram inside the plane. The theorem also holds 
for oriented links and diagrams. One then has to take into account
all possible coherent orientations of the diagrams involved in the moves.
\end{theorem}

\ \\
\centerline{\psfig{figure=R1R2R3Trieste.eps,height=11.5cm}}\ \\
\centerline{Figure 1.6; Three Reidemeister moves: $R_1,R_2$ and $R_3$}\ \\

\subsection{Fox 3-colorings of link diagrams}\ \\

The simplest invariant of links which distinguishes the trefoil knot and the trivial knot
({\parbox{2.1cm}{\psfig{figure=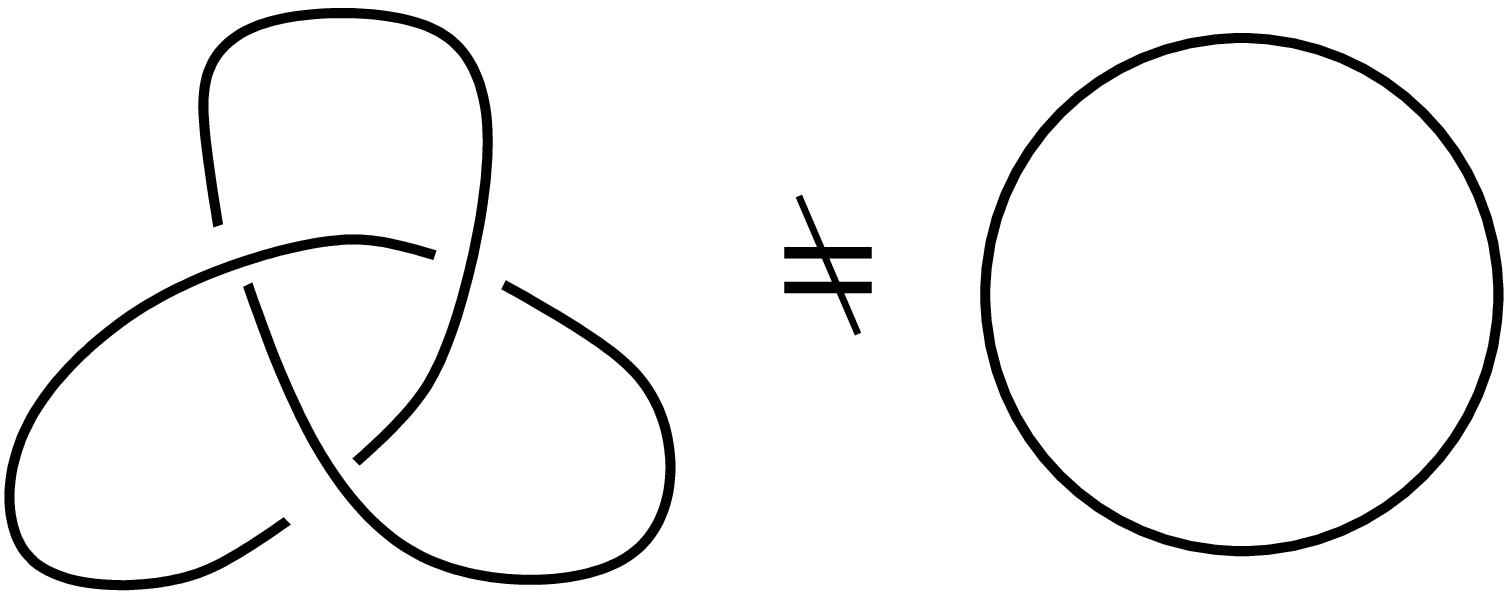,height=0.8cm}}}) is 
the Fox tricoloring invariant (denoted  $tri(L)$).  
It is an invariant which does not require much more than counting.
The idea of tricoloring was introduced by Ralph Hartzler Fox (1913 -1973)
around 1956 when he was
explaining knot theory to undergraduate students
at Haverford College (``in an attempt to make the subject
accessible to everyone" \cite{C-F}); [C-F,Chapter VI,Exercises 6-7], \cite{Fo-2}.
 It was also popularized in articles directed toward 
 middle and high school teachers  and students \cite{Cr,Vi,P-6}.

\begin{definition}[{\cite{P-1}}]\label{Definition 1.4}
We say a link diagram $D$ is Fox tricolored if every arc is colored $r$
(red), $b$ (blue) or $y$ (yellow) ( we consider arcs of the diagram literally,
so that in the undercrossing one arc ends and the second starts;
compare Fig.1.7, 1.9), and at any given crossing either all three colors
appear or only one color appears.
The number of different Fox tricolorings is denoted by $tri(D)$. If a tricoloring
uses only one color we say that it is a trivial Fox tricoloring.
\end{definition}

\centerline{\psfig{figure=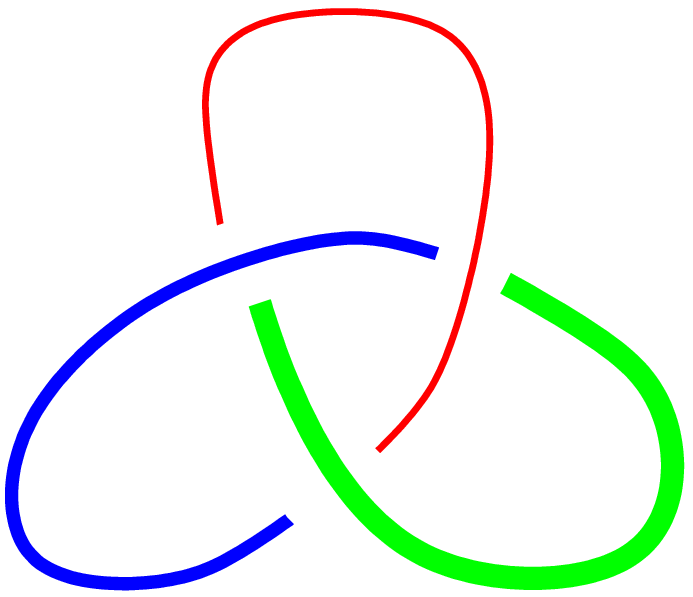,height=4.1cm}}
\begin{center}
Fig. 1.7. Different colors are marked by lines of different thickness.
\end{center}

\begin{proposition}\label{Proposition 1.7}
The number of Fox tricolorings of $D$, $tri(D)$ is an (ambient isotopy) link invariant. 
In particular, the tricolorability, that 
is  the existence of a non-trivial Fox tricoloring, is a link invariant.
\end{proposition}

\begin{proof}

We have to check that $tri(D)$ is preserved under the Reidemeister moves.
The invariance under $R_1$ and $R_2$ is illustrated in Fig.~1.8, and the
invariance under $R_3$ is illustrated in Fig.~1.9.
\end{proof}
\ \\
\centerline{\psfig{figure=R1R23-col.eps}}
\centerline{Fig. 1.8}

\ \\

\centerline{\psfig{figure=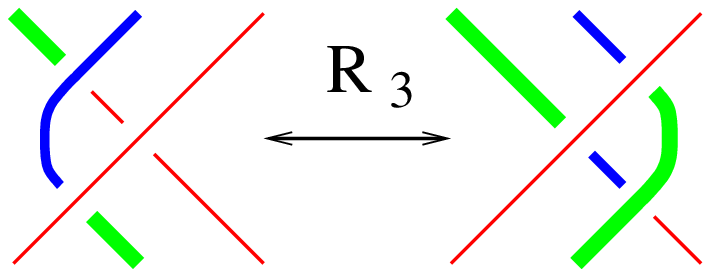}}
\begin{center}
Fig. 1.9
\end{center}

Because the trivial knot has only trivial tricolorings, $tri(T_1)=3$,
and the trefoil knot allows a nontrivial tricoloring (Fig.1.7),
it follows that the trefoil knot is a nontrivial knot.

\begin{exercise}\label{1.3}
Find the number of tricolorings for the trefoil knot ($3_1$),  the figure
eight knot\footnote{The figure eight knot is often called the Listing knot, as Listing noticed in 1849 that 
it is equivalent to its mirror image. The notation $4_1$ refers to the fact that it is the first knot 
of 4 crossings in knot tables.} ($4_1$), and the square knot ($3_1\# \bar 3_1$, see Fig.1.10).
Then deduce that these knots are pairwise different.
\end{exercise}

It is very difficult to prove any nontrivial result using our previous definition of tricoloring. 
For example how would you prove the following statement?
\begin{proposition}\label{Proposition 1.7 Trieste}
$tri(L)$ is always a power of $3$.
\end{proposition}
We can see immediately that if we tricolor arcs of a diagram $D$ without Fox coloring conditions we 
get $3^{\lambda}$ possibilities, where $\lambda$ is the number of arcs of $D$. 
Thus for a diagram without a crossing proposition \ref{Proposition 1.7 Trieste} 
holds but if $D$ has a crossing we only can say that $tri(D)\leq 3^{\lambda}$

Proposition \ref{Proposition 1.7 Trieste} becomes easy to prove if we introduce some basic  language 
of linear algebra or 
abstract algebra. Namely:\\
\begin{proof}  
Denote the colors of the Fox tricoloring by $0,1$ and $2$ and treat them modulo $3$, that is, 
as elements of the group (or field) $\Z_3$. All colorings of the arcs of
a diagram using colors $0,1$,and $2$ (not necessarily permissible Fox tricolorings)
can be identified with the group $\Z_3^{\lambda}$ (or the linear space over $\Z_3$).
The (permissible) Fox tricolorings can be characterized by the property that at
each crossing, the sum of the colors is equal to zero modulo $3$.
Thus Fox tricolorings form a subgroup (linear subspace) of $Z_3^{\lambda}$. We denote this group $Tri(D)$. 

I encourage you to play around with this concept. Notice that 
trivial colorings form a one dimensional subspace, so one can should consider the quotient space of 
all Fox 3-colorings by the subspace of trivial tricolorings $\Z_3^{tr}$. 
We call this quotient space the space of reduced Fox tricolorings; 
$Tri^{rd}(D)= Tri(D)/\Z^{tr}_3$.
\end{proof}
Given our  an easy success with the proof of Proposition  \ref{Proposition 1.7 Trieste}, let us 
try our skills on the following fact and its useful corollary. 
Recall that an $n$-tangle is a part of a link diagram placed in a 2-disk with
$2n$ points on the disk boundary: $n$ inputs and $n$ outputs (however only if a tangle is oriented we 
have unique notion of inputs and outputs); see examples in Figures 1.10 -- 1.12. 

\begin{proposition}\label{Proposition 1.8 Trieste}
\begin{enumerate}
\item[(i)] For any Fox 3-coloring of a 1-tangle; see Fig. 1.12(a), boundary arcs share a color .
\item[(ii)] $tri(L_1)tri(L_2)=3tri(L_1\# L_2)$, where $\#$ denotes
the connected sum of links\footnote{A diagram $D_1\# D_2$ is a connected sum of diagrams $D_1$ and $D_2$ if there 
is a simple closed curve cutting $D_1\# D_2$ in exactly two points and $1$-tangles obtained by 
cutting $D_1\# D_2$ by the curve have $D_1$ and $D_2$ as their closures. A link $L_1\# L_2$ is 
a connected sum of links $L_1$ and $L_2$ if there is a diagram of $L_1\# L_2$ which is a 
connected sum of diagrams of $L_1$ an $L_2$. Connected some maybe not unique, and may 
depend on components of links connected in the connected sum and on orientation of links.}; see Fig. 1.10.
\end{enumerate}
\end{proposition}
\begin{proof}
(i) Let $T$ be our Fox tricolored tangle and let the 1-tangle $T'$
be obtained from
$T$ by adding a trivial component $C$ below $T$, close enough to the boundary
of the tangle, so that it cuts $T$ only near the boundary points; Fig.1.11(b).
Obviously the tricoloring of $T$ can be extended to a tricoloring of $T'$ (in
three different ways) because
the tangle $T'$ is ambient isotopic to a tangle obtained from $T$ by adding
a small trivial
component disjoint from $T$. However, if we try to color $C$,
we see immediately that it
is possible if and only if the input and the output arcs of $T$ have the same color. Namely, if 
$x$ is the color of a point on $C$ and $a$ and $b$ colors of the input and the output then following $C$ 
and using Fox tricoloring rules at two crossings of $C$ with $T$ we get $x= a-b+x$, so $a=b$; see Figure 1.11.\\
(ii) If we consider the connected sum $L_1\# L_2$, we see from the part (i) that the arcs joining $L_1$ and $L_2$     
have the same color. Therefore the formula 
$tri(L_1\# L_2)=\frac{1}{3}tri(L_1)tri(L_2)$ follows.
\end{proof}

\ \\
\ \\
\centerline{\psfig{figure=connectedsum.eps,height=2.6cm}}
\begin{center}
Fig. 1.10; connected sum of link diagrams
\end{center}
\ \\
\centerline{\psfig{figure=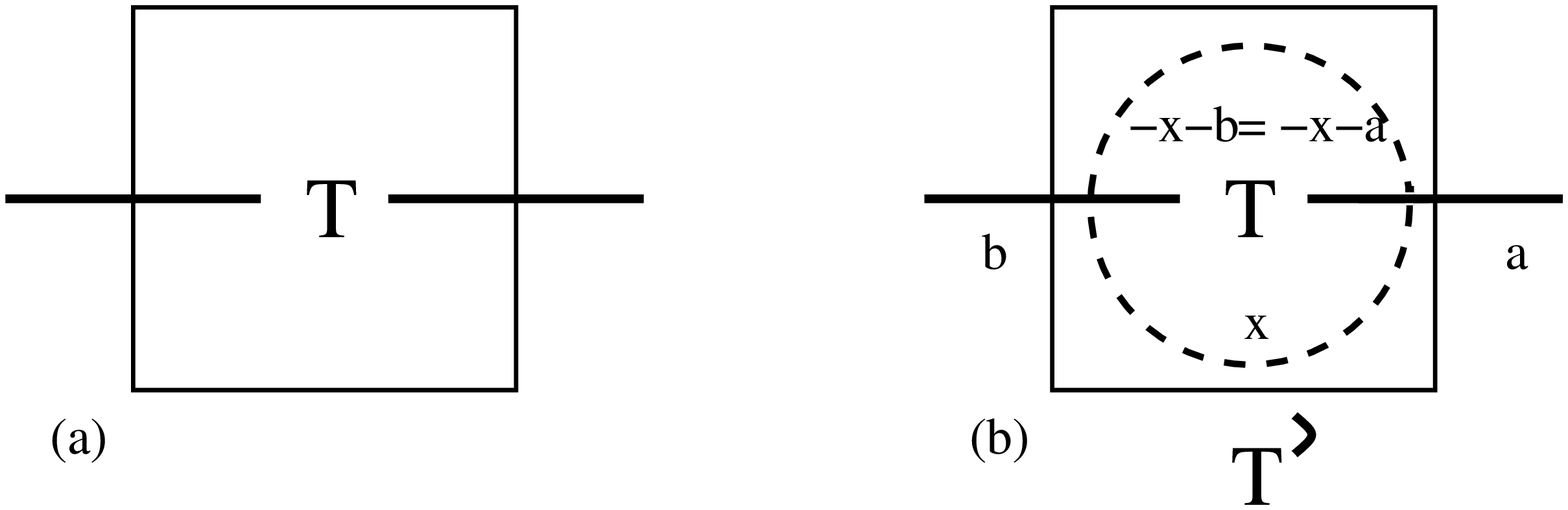,height=3.7cm}}
\begin{center} Fig. 1.11; 1-tangle and tricoloring of its boundary points; $a=b$  \end{center}

\ \\
The next proposition gives a very interesting property relating the number of Fox tricolorings  
of four unoriented links which differ in a small neighborhood as in Figure 1.12. 
Using basic algebra we can only partially prove Proposition \ref{Proposition 1.9 Trieste}. 
Tomorrow I will show you how to place more structure 
on colorings (symplectic structure) to fully prove the proposition and its generalizations. 

\begin{proposition}\cite{P-3}\label{Proposition 1.9 Trieste}
Let $D_+, D_-, D_0$ and $D_{\infty}$ denote four unoriented link diagrams
(of links $L_+, L_-, L_0$ and $L_{\infty}$) as in Fig. 1.12. Then among the four numbers $tri(L_+), tri(L_-), tri(L_0)$ and
$tri(L_{\infty})$, three are equal and the fourth is
 $3$ times bigger then the rest.
\end{proposition}

\centerline{\psfig{figure=L+L-L0Linf.eps}}
\begin{center}
Fig. 1.12 
\end{center}

We first prove here the weaker fact that among these four numbers\\
 either all 4 are  equal or 3 of them are equal and the 4'th is $3$ times bigger then the rest
(the rest of the proof will wait till tomorrow). 

\begin{proof}
Consider a crossing $p$ of the diagram $D$. If we cut  a neighborhood of $p$ out of $D$,
 we are left with the 2-tangle $T_D$ (see Fig.1.13(a)).
The set of Fox tricolorings of $T_D$, $Tri(T_D)$, forms a linear space over $\Z_3$ with 
subspaces $Tri(D_+),Tri(D_-), Tri(D_0)$ and $Tri(D_{\infty})$.
Let $x_1,x_2,x_3,x_4$
be elements of $Tri(T_D)$ corresponding to arcs cutting the boundary of
the tangle; see Fig.1.13(b).
Then any element of $Tri(T_D)$ satisfies the equality \ $x_1-x_2+x_3-x_4=0$.
To show this, we proceed as in part (i) of Proposition \ref{Proposition 1.8 Trieste},
 see Figure 1.13(b).
Any  element of $Tri(D_+)$ (resp. $Tri(D_-)$,
$Tri(D_0)$ and $Tri(D_{\infty})$) satisfies additionally the equation \
$x_1=x_3$ (resp. $x_2=x_4$, $x_2=x_3$, and $x_2=x_1$). Thus $Tri(D_+)$ (resp.
$Tri(D_-)$, $Tri(D_0)$ and $Tri(D_{\infty})$) is a subspace of $Tri(T_D)$ of
codimension at most one.
Let $F$ be the subspace of $Tri(T_D)$ given by the equations
$x_1=x_2=x_3=x_4$, that is, the space of 3-colorings monochromatic on
the boundary of
the tangle. $F$ is a subspace of codimension at most one in any
of the spaces
$Tri(D_+)$, $Tri(D_-)$, $Tri(D_0)$, $Tri(D_{\infty})$. Furthermore
the common part of any two of $Tri(D_+)$, $Tri(D_-),Tri(D_0), Tri(D_{\infty})$
is equal to $F$. To see this, we just compare the defining relations
for these spaces.
Finally, notice that
$Tri(D_+)\cup Tri(D_-)\cup Tri(D_0)\cup Tri(D_{\infty})=Tri(T_D)$.

We have the following possibilities:
\begin{enumerate}
\item
[(1)] $F$ has codimension 1 in $Tri(T_D)$.\\
Then by the above considerations:\\
One of $Tri(D_+),Tri(D_-),Tri(D_0), Tri(D_{\infty})$
is equal to $Tri(T_D)$. The  remaining three spaces are equal to $F$ and 
 Proposition \ref{Proposition 1.9 Trieste} holds.
\item
[(2)] $F=Tri(D_+)=Tri(D_-)=Tri(D_0)= Tri(D_{\infty})=Tri(T_D)$,
\item
[(3)] $F$ has codimension $2$ in $Tri(T_D)$. Then
$3|F|=tri(D_+)=tri(D_-)=tri(D_0)= tri(D_{\infty})=\frac{1}{3}tri(T_D)$
\end{enumerate}
This completes the weaker statement of Proposition \ref{Proposition 1.9 Trieste}. 
To prove Proposition \ref{Proposition 1.9 Trieste} fully, one must exclude cases (2) and (3).
To exclude (2) and (3) one can use the Goeritz matrix of the link diagram; see \cite{P-Book}.
In the second part we show how to use the concept of Lagrangian tangles to show essential 
generalization of Proposition \ref{Proposition 1.9 Trieste} (the concept was introduced in \cite{DJP}).

\end{proof}
\centerline{\psfig{figure=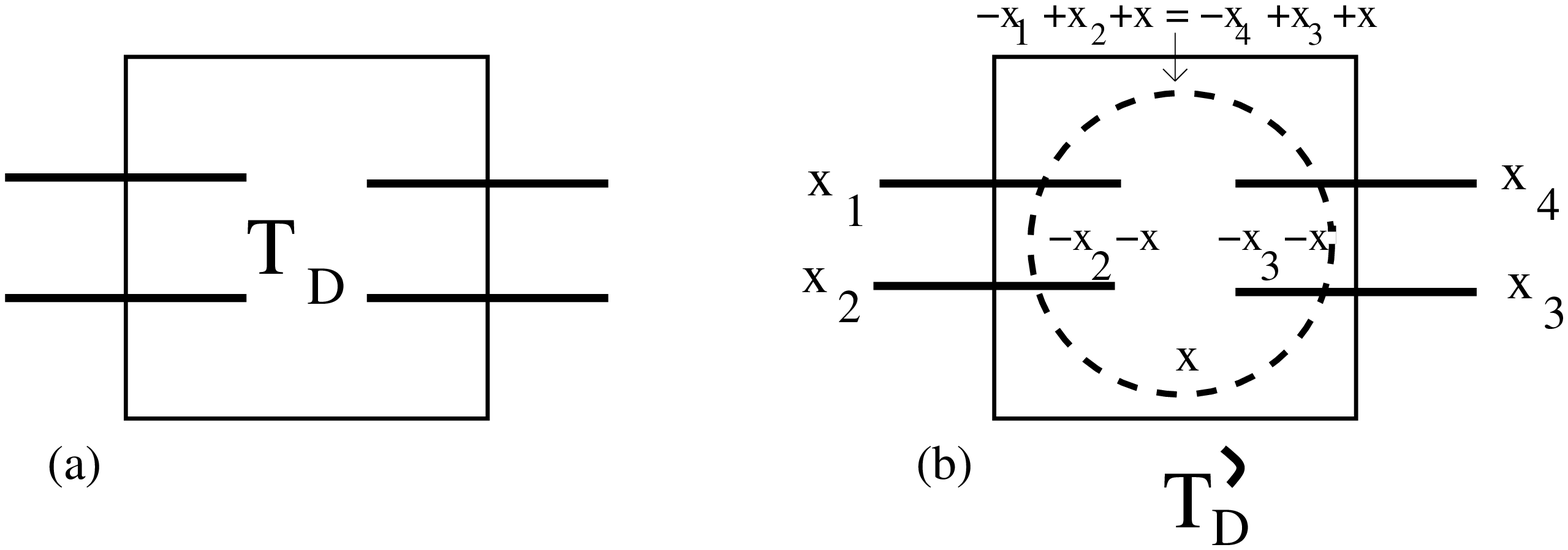,height=4.6cm}} 
\begin{center} Fig. 1.13; 1-tangle and tricoloring of its boundary points; $x_1-x_2+x_3-x_4=0$  \end{center}

We will show below that the dimension of the space of Fox 3-colorings of a link is bounded from above 
by the bridge index of the link. For this we need few basic definitions: \\
Let $L$ be a link embedded in $R^3$ which meets a plane $E\subset R^3$ in $2k$
points such that the arcs of $L$ contained in each half-space relative to $E$
possess orthogonal projections onto $E$ which are simple and disjoint.
$(L,E)$ is called a $k$-bridge presentation of $L$; \cite{B-Z}. The bridge index of a link $L$, denoted $bridge(L)$,
 is a minimal number $k$ such that $L$ has a $k$-bridge presentation. 
Notice that the $k$-bridge presentation of $L$ 
can be interpreted as an embedding of $L$ with exactly $k$ minima and $k$ maxima (in the $z$ direction).
\begin{proposition}\label{Proposition 1.10 Trieste}
For any link $L$ we have $$tri(L)\leq 3^{bridge(L)}$$.
\end{proposition}
\begin{proof} If we color the bridges of a diagram, then the 3-coloring of the other arcs is 
 uniquely determined. 
It may happen however, that we get ``contradictions" at some minima; which leads to 
 the inequality in Proposition \ref{Proposition 1.10 Trieste}.
\end{proof}

\begin{remark}\label{Remark 1.11 Trieste}
We can look at links or tangles with $n$ bridges from a different perspective, by 
organizing diagrams along the $y$ axis that, is we deal with maxima (and  minima) along the $y$ axis. 
In the case of a $2$-tangle we also have $4$ minimal (boundary) points, in addition 
to $n$ maxima ($\cap$) and $n-2$ minima ($\cup$); compare Figure 1.14.
\\ \ \\ \
\centerline{\psfig{figure=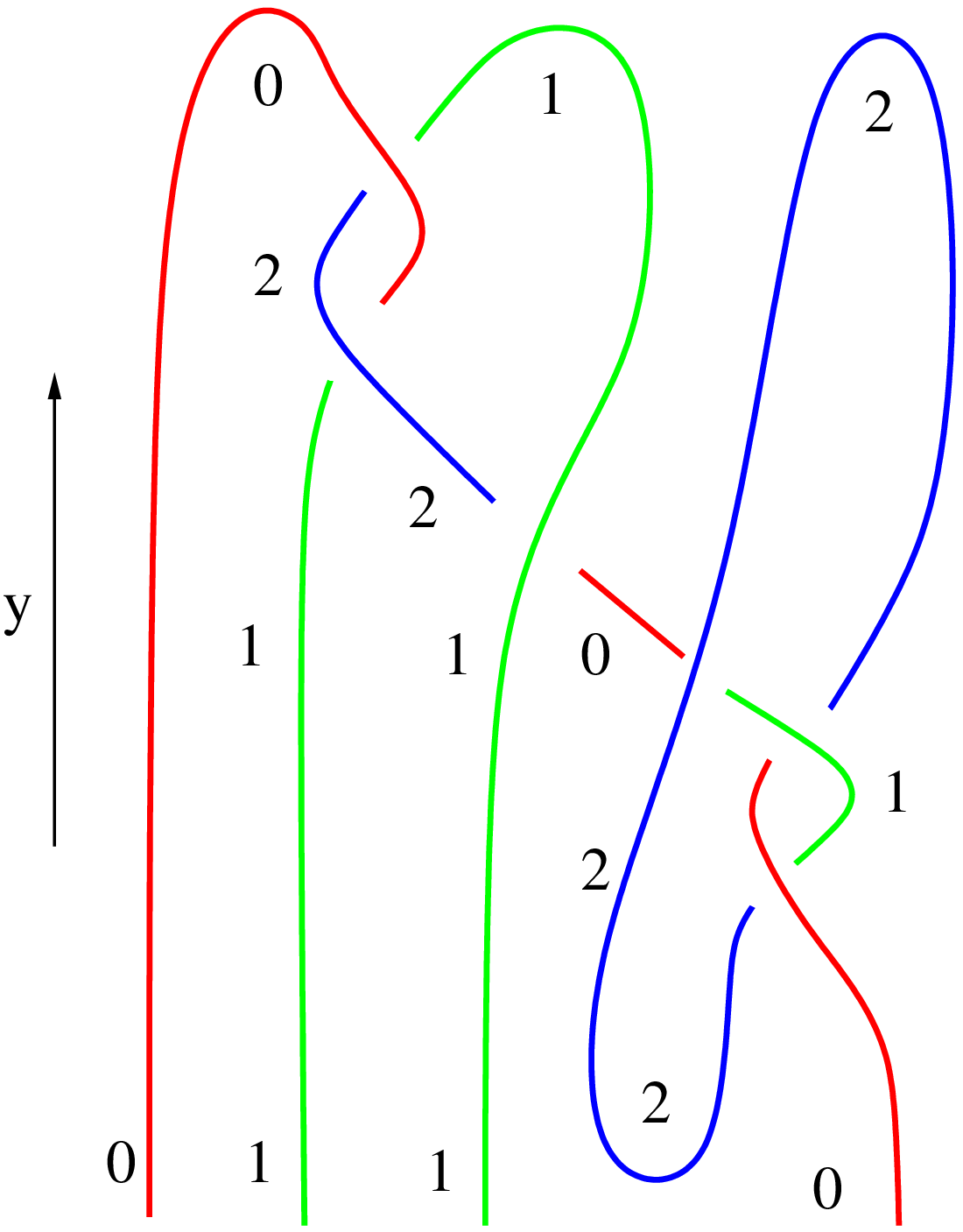,height=8.5cm}} 
\begin{center} Fig. 1.14; Diagram of a 2-tangle with 3 maxima \end{center}
 
We observe that that if we 
tricolor maxima, it will propagate until we reach minima ($\cup$) which will give obstruction (additional 
 relations) to possible 3-colorings. If we start with $n$-maxima, we also start with $\Z_3^n$ as the space 
of colorings. When we move along our diagram down, with respect to the $y$ axis  (like a braid) we uniquely 
color the arcs of the diagram and, at any level, keeping n-dimensional space of boundary colorings
 until we reach minima.
Assume we deal with a 4-tangle $T_D$. Then we have $n-2$ minima leading to $n-2$ relations on 
 $\Z_3^n$. Thus we are left with at least a 2-dimensional space. In Section 2 we show more: 
  the colorings of boundary points span exactly a 2-dimensional space. 
To express this algebraically, we consider a linear map $\psi: Tri(T_D) \to \Z_3^4$, in 
which the coloring of the 2-tangle $T_D$ yields a coloring of the four boundary points. If we start from 
an $n$-tangle with $n$-maxima, we have an isomorphism $Tri(T_D) \to \Z_3^n$ and after adding 
the $n-2$ relations, the image of $\psi$ is 2-dimensional.

In conclusion, this shows that any 2-tangle has a 3-coloring that is not monochromatic on the boundary. 
This will be discussed, given additional tools, more generally, tomorrow\footnote{Tomorrow's 
talk will introduce a symplectic structure on the space of colorings of a tangle boundary  
which does not apply to virtual links and tangles (which we heard about today). Therefore one should  
mention that the considerations in the observation above apply partially to virtual tangles as well. 
On the other hand, part of the proof of Proposition \ref{Proposition 1.9 Trieste} does not work
for virtual links: the equality 
$x_1-x_2+x_3-x_4=0$  does not always hold for diagrams with virtual crossings, and is related with 
the fact that a virtual crossing alone does not satisfies the property for any nontrivial coloring. 
For the virtual crossing 
{\parbox{1.2cm}{\psfig{figure=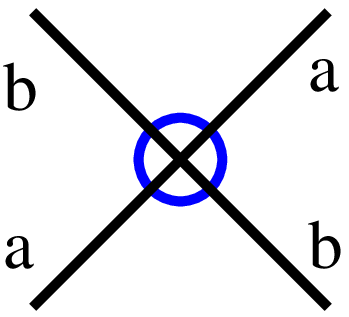,height=0.8cm}}} we have $a-b+a-b=2(a-b)$. In the 
virtual knot theory we have two forbidden moves as an arc cannot be moved under or over 
  a virtual crossing (the first forbidden move:  
{\parbox{2.2cm}{\psfig{figure=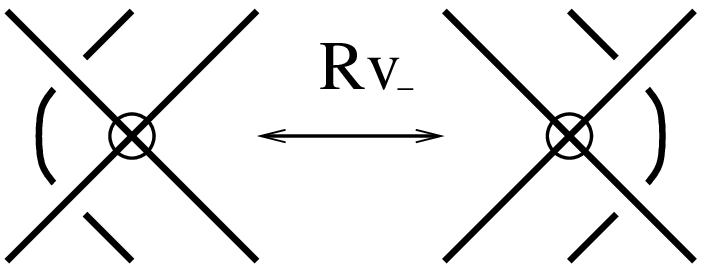,height=0.8cm}}}) and the second forbidden move:
{\parbox{2.2cm}{\psfig{figure=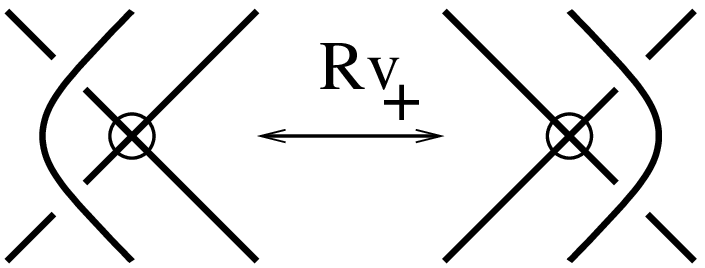,height=0.8cm}}).
The theory of Fox colorings works for virtual links or tangles. It works also if the second 
forbidden move is allowed (so can be used for welded knots 
described in Kauffman's talk). Fox colorings and more generally quandle colorings (see Section 3)
are preserved by this move.}}.
\end{remark}
\subsection{Fox 3-colorings and the Jones polynomial}\ 

In many talks we heard about the Jones polynomial -- the great breakthrough in knot theory, in 1984.
 
I noticed a connection between Fox tricolorings and the Jones polynomial
when I analyzed the influence of 3-moves on 3-colorings and
the Jones polynomial \cite{P-1}.
\begin{definition}\label{Definition 1.8}\ \\
The local change in a link diagram which replaces parallel lines
by $n$ positive half-twists is called an n-move; see Fig.1.14.
\end{definition}

\ \\
\centerline{\psfig{figure=n-moves.eps}}
\begin{center}
Fig. 1.14; 3-move and n-move
\end{center}

\begin{lemma}\label{Lemma 1.13 Trieste}
 Let the diagram $D_{+++}$ be obtained from $D$ by a 3-move
(Fig.1.14(a)). Then:
\begin{enumerate}
\item
[(a)] $tri(D_{+++})=tri(D)$,
\item
[(b)] $V_{D_{+++}}(e^{2\pi i/6})=\pm i^{(com(D_{+++})-com(D))}V_D(e^{2\pi i/6})$, where $V$ is the Jones polynomial, 
and $com(D)$ denotes the number of link components of $D$,
\item
[(c)] $F_{D_{+++}}(1,-1)=  F_D(1,-1)$,  where $F$ is the Kauffman polynomial.
\end{enumerate}
\end{lemma}

Before we prove Lemma \ref{Lemma 1.13 Trieste} 
let us recall definition of the Jones polynomial (1984) and the specialization 
of the Kauffman polynomial first introduced by Brandt-Lickorish-Millett and Ho \cite{BLM,Ho} (1985).

\begin{definition}
 \begin{enumerate}
 \item[(J)] The Jones polynomial
$V_L(t)$ of an oriented link $L$ is a link invariant ($V_L(t)\in \Z[t^{\pm 1/2}]$) normalized to be one for
the trivial knot and satisfies the skein relation
$$ t^{-1}V_{{\psfig{figure=L+maly.eps}}}(t) -
tV_{{\psfig{figure=L-maly.eps}}}( t) =
(t^{\frac{1}{2}} - t^{-\frac{1}{2}})
V_{{\psfig{figure=L0maly.eps}}}(
t).$$
 \item[(K)] The Brandt-Lickorish-Millett-Ho polynomial $Q_L(x)$ is normalized to be one for
the trivial knot and satisfies the skein relation
$$Q_{L_{{\psfig{figure=L+nmaly.eps,height=0.4cm}}}}(x) +
Q_{L_{{\psfig{figure=L-nmaly.eps,height=0.4cm}}}}(x) =
xQ_{L_{{\psfig{figure=L0nmaly.eps,height=0.4cm}}}}(x) +
xQ_{L_{{\psfig{figure=Linftynmaly.eps,height=0.4cm}}}}(x).$$
The Kauffman 2-variable polynomial $F(a,x)$ satisfies $F(1,x)= Q_L(x)$.
\end{enumerate}
\end{definition}
\begin{exercise}
\begin{enumerate}
\item[(i)] Show that $V_{T_n}=(-t^{1/2} - t^{ -1/2})^{n-1}$, for $T_n$ being the trivial link 
of $n$ components.
\item[(ii)] Show that $V_L(t) \in \Z[t^{\pm 1}]$ if $L$ has odd number of components and 
$t^{1/2}V_L(t) \in \Z[t^{\pm 1}]$ if $L$ has even number of components.
\item[(iii)] Show that $V_K(t) - 1$ is divisible by $(t-1)(t^3-1)$ for any knot $K$.
\item[(iv)] Show that $V_L(t) - V_{T_{com (L)}}$ is divisible by $(t^3-1)$ for any link $L$.
Here $com (L)$ denotes the number of components of $L$ and $T_k$ is the trivial link of $k$ components 
\end{enumerate}
\end{exercise}
{\bf Proof of Lemma \ref{Lemma 1.13 Trieste}}.\\
We prove (a) and (c) and partially (b) (one of two possible orientation choices).
\begin{enumerate}
\item
[(a)]
The bijection between 3-colorings of $D$ and $D_{+++}$ is illustrated in
Fig. 1.15.
\ \\ \ \\ \ \\
\centerline{\psfig{figure=3-movecol.eps}}
\begin{center}
Fig. 1.15
\end{center}
\item
[(c)] $F_{D_{+++}}(1,-1)=-F_{D_{+}}(1,-1)-F_{D_{++}}(1,-1)-
F_{D_{\infty}}(1,-1) =
-F_{D_{+}}(1,-1)+ F_D(1,-1)+F_{D_{+}}(1,-1)+F_{D_{\infty}}(1,-1)-
F_{D_{\infty}}(1,-1)= F_D(1,-1)$.
\item
[(b)] Assume that arcs in Figure 1.15(a) have parallel orientation. Then for $t=e^{2\pi i/6}$ 
($t^{1/2}= e^{\pi i/6}$) we have:\\
$V_{D_{+++}}= t^2V_{D_{+}} + t(t^{1/2}-t^{-1/2})V_{D_{++}}= t^2V_{D_{+}} + t^3(t^{1/2}-t^{-1/2})V_D + 
t^2(t^{1/2}-t^{-1/2})^2V_{D_{+}}= t^2(t-1 +\frac{1}{t})V_{D_{+}} + t^{1/2}(t^3-t^2)V_D = 
 \frac{t^3+1}{t+1}V_{D_{+}} + t^{3/2}\frac{t^3+1}{t+1}V_{D} - t^{3/2}V_{D} = -e^{\pi i/2}V_{D} = 
-iV_{D} $, as needed.\\
In the case when a 3-move is not preserving orientation, we would have to consider several involved 
cases, but we can make shortcut using so called Jones reversing result\footnote{It was initially 
proven in a series of involved papers but now it has an easy proof using 
the Kauffman bracket polynomial which do not depend on a link orientation; compare \cite{P-Book}. 
Precisely, we have: Suppose that $L_i$ is a component of an oriented
link $L$ and $\lambda = \lk (L_i , L-L_i)$. If $L'$ is a link obtained from $L$ by reversing
the orientation of the component $L_i$ then $V_{L'}(t) = t^{-3\lambda}V_L(t)$.}, 
that is if one changes 
an orientation of some components of a link, then its Jones polynomial is changed in a precisely 
described way, in particular by multiplying by a number being the power of $t^{3}$ 
(in our case the power of $-1$).

\end{enumerate}

One can easily check that for a trivial $n$-component
link, $T_n$,
$tri(T_n)=3^n=3V_{T_n}^2(e^{2\pi i/6})=3(-1)^{n-1}F_{T_n}(1,-1)$. Furthermore
it follows from Lemma 1.9 that as long as a link $L$ can be obtained from a trivial link
by 3-moves we have:
$tri(L)=3|V_L^2(e^{2\pi i/6})|=3|F_L(1,-1)|$.

These immediately lead to three questions:
\begin{enumerate}
\item[(1)] ({\bf Montesinos-Nakanishi 3-move conjecture}). \\
Any link can be reduced to a trivial link by a finite sequence of 3-moves.
\item[(2)] Is it true that $tri(L) = 3|F_L(1,-1)|$?
\item[(3)] Is it true that  $tri(L)=3|V_L^2(e^{2\pi i/6})|$?
\end{enumerate} 

Formulas (2) and (3) follow immediately from (1) as equalities from (2) and (3) hold for 
trivial links and it is propagated by 3-moves. Thus we proved (2) and (3) for any link 
which can be reduced by 3-moves to a trivial link.

Y.Nakanishi first considered the conjecture (i) in 1981. J.Montesinos
analyzed 3-moves before, in connection with 3-fold dihedral branch
coverings, and asked a related but different question.
The conjecture was proved in many special cases (e.g. \cite{Che}) but it was an 
open problem for over 20 years. In 2002 it was showed by M.K.D{\c a}bkowski and the author
that the conjecture does not hold. The smallest counterexample we found, suggested first by Q.~Chen, 
has 20 crossings, see Figure 1.16, \cite{D-P-1}.. We conjecture that it is in fact the smallest counterexample,
that is every link up to 19 crossings can be reduced to a trivial link by 3-moves, furthermore we 
predict that every link of 20 crossing is reduced by a 3-move either to the trivial link or to the 
Chen link (up to the mirror image).  With todays computers it should be laborious but doable exercise -- please try it!
\\ \ \\
\centerline{\psfig{figure=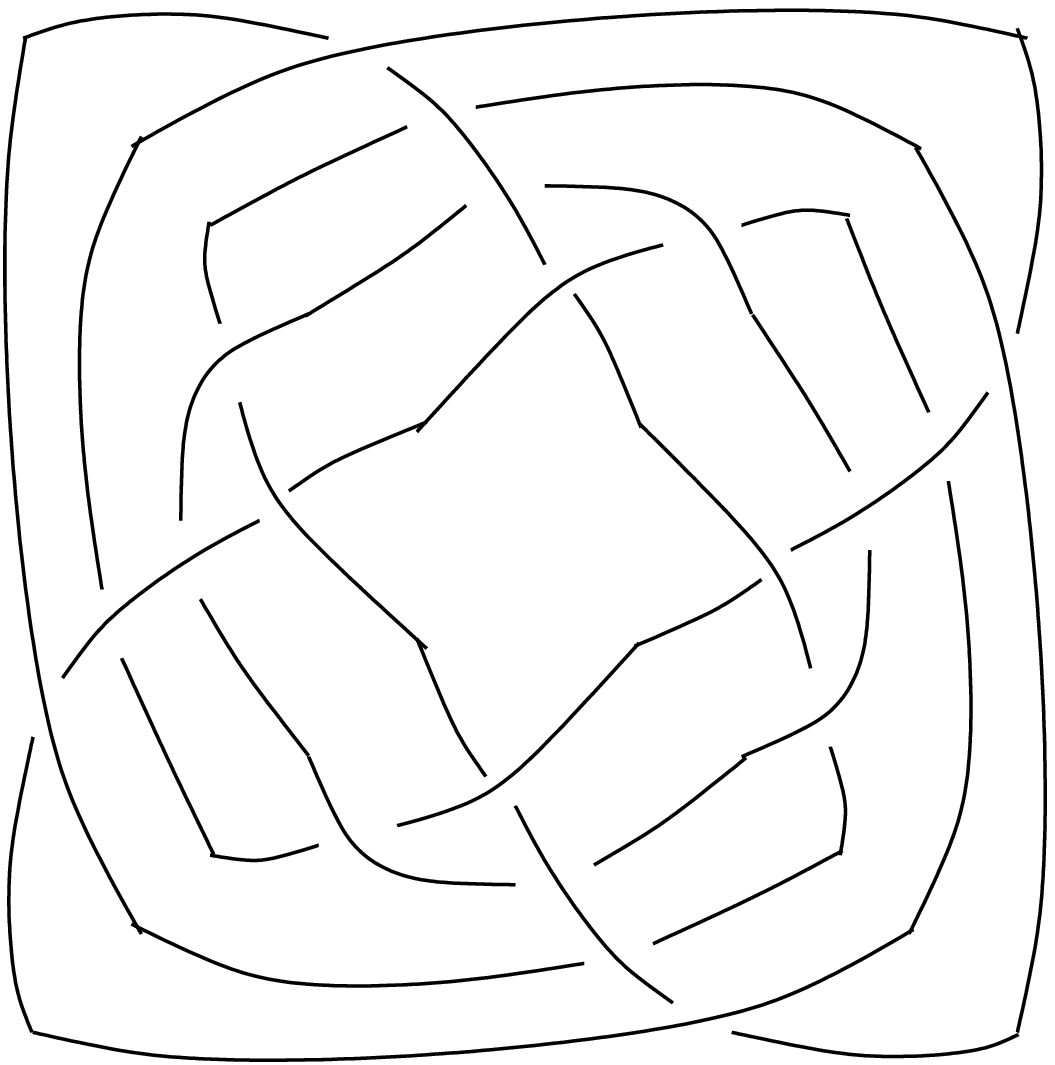,height=4.2cm}}
\centerline{Fig. 1.16; the Chen link, the closure of the 5-string 
braid $(\sigma_2\sigma_1^{-1}\sigma_2\sigma_3\sigma_4^{-1})^4$}
\\ \

The Montesinos-Nakanishi conjecture does not hold but  the formulas (2) and (3) linking 
tricoloring with the Jones and Kauffman polynomials holds for any link. 

The proof of (a) in \cite{P-1} uses Fox's interpretation of 3-coloring
and the connection with the first homology group of the
branched 2-fold cover of $S^3$
branched over the link. However, a simple, totally elementary proof follows from Proposition \ref{Proposition 1.9 Trieste}.\ \\
\begin{proof}
Because $tri(L)$ is a power of $3$, we can consider the signed version of the
tricoloring defined by:\ \  $tri'(L)=(-1)^{log_3(tri(L))}tri(L)$.\
It follows from Proposition \ref{Proposition 1.9 Trieste} that $$tri'(L_+)+tri'(L_-)=
-tri'(L_0)-tri'(L_{\infty}).$$
This is however exactly the recursive formula for the Kauffman polynomial
$F_L(a,x)$ at $(a,x)=(1,-1)$.
Comparing the initial data (for the unknot) of $tri'$ and $F(1,-1)$
we get generally that:  \ $-3F_L(1,-1)=tri'(L)=(-1)^{log_3(tri(L))}tri(L)$,
which proves part (b) of Theorem 1.13.
Part (a) follows from Lickorish's observation \cite{Li}, that
$F_L(1,-1)=(-1)^{com(L)}V_L^2(e^{2\pi i/6})$.
This observation can be directly proven from the Kauffman bracket polynomial 
version of the Jones polynomial. For people who attended Lou Kauffman talk it should 
be a pleasure exercise: just consider the difference of squares of the Kauffman bracket relation for 
$L+$ and $L_-$ (that is $\langle$ {\parbox{0.5cm}{\psfig{figure=L+nmaly.eps}}}
$\rangle^2 - \langle$ {\parbox{0.5cm}{\psfig{figure=L-nmaly.eps}}}
$\rangle^2$). You will get the relation of the, so called, Dubrovnik version of the Kauffman polynomial which 
can be converted to the standard one.
\end{proof}

I would challenge you to find completely elementary proof of Proposition \ref{Proposition 1.9 Trieste} 
or directly formulas (b) and (c) (as we noted all three facts are related by elementary consideration). 
As a prize I offer a copy of my book \cite{P-Book}. 


Tomorrow I will define general Fox k-colorings and Fox coloring group, and 
I will place the theory of Fox coloring in more
general (sophisticated) context, and apply it to the analysis
of k-moves (and rational and braid moves) of $n$-tangles.
Interpretation of tangle colorings as Lagrangians in symplectic
spaces is our main (and new) tool. In the second lecture tomorrow,
I will also mention another motivation for studying 3-moves: to understand
skein modules based on their deformation.

\section{ Fox colorings, rational moves, and Lagrangian tangles}

Many of you, likely, wondered yesterday why we consider only  3-colorings not, say 
generally $n$-colorings. Some of you probably tried to replace the relation 
$a+b+c \equiv 0 \mod 3$ by the relation $a+b+c \equiv 0 \mod k$, and noticed that it does not 
work well with Reidemeister moves. In fact, as observed by Fox, the proper relation to 
generalize is $2b-a-c \equiv 0 \mod 3$. This leads to Fox $k$-colorings:
\begin{definition}\label{2.1}
\begin{enumerate}
\item [(i)]
We say that a link (or a tangle) diagram is Fox k-colored if every
arc is colored by one of the numbers $0,1,...,k-1$ (forming a
group $Z_k$) in such a way that
at each crossing the sum of the colors of the undercrossings is equal
to twice the color of the overcrossing modulo $k$;algebraically $c\equiv 2b-a \mod k$ as illustrated in Fig.2.1.
\item [(ii)]
The set of Fox $k$-colorings forms an abelian group (or $\Z_k$-module), denoted by $Col_k(D)$.
The cardinality of the group will be denoted by $col_k(D)$.
For an $n$-tangle $T$ each Fox $k$-coloring of $T$ yields a
coloring of boundary points of $T$ and we have the homomorphism
$\psi :Col_k(T) \rightarrow Z_k^{2n}$
\end{enumerate}
\end{definition}.
\centerline{\psfig{figure=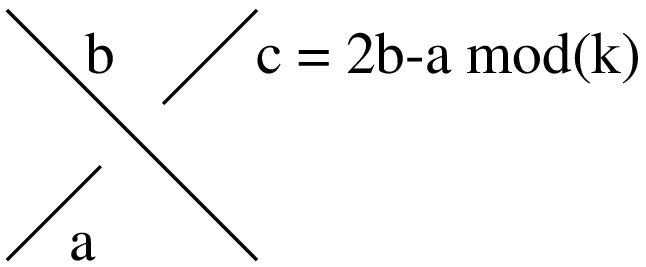,height=2.6cm}}
\centerline{Fig. 2.1}
It is a pleasant exercise to show that $Col_k(D)$ is unchanged
by Reidemeister moves (see Figure 2.2),

\ \\
\centerline{\psfig{figure=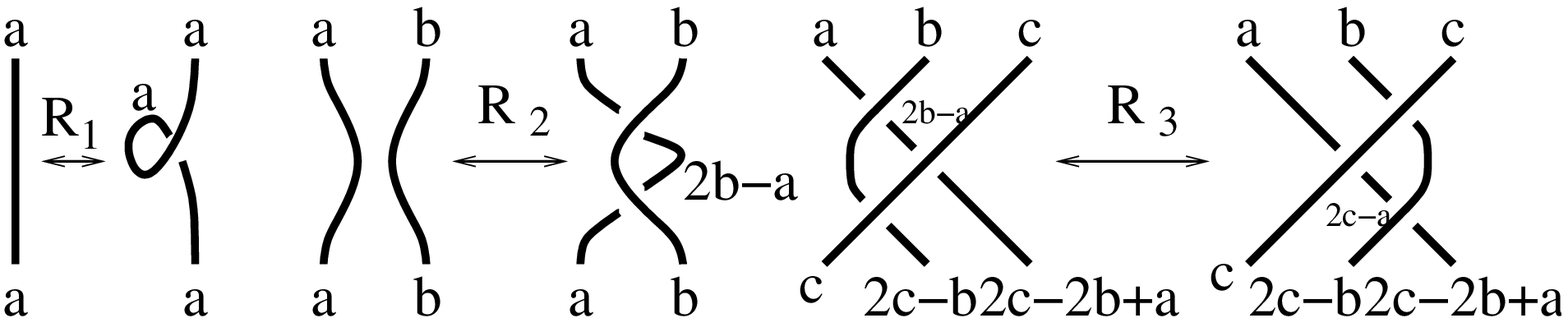,height=2.6cm}}
\centerline{Fig. 2.2}
\\

I will start this part from the basic observations on Fox $k$-colorings 
analogous to those proven yesterday for 3-colorings.
 The talk will culminate by the
introduction of the symplectic structure on the boundary
of a tangle in such a way that tangles yields Lagrangians
in the symplectic space. We end with some corollaries, in particular the method 
to recognize often that a virtual tangle is not a classical tangle (by boundary $k$-coloring  
comparison). 

We follow here \cite{P-3} and \cite{DJP} (see \cite{P-5} for historical introduction).

\begin{proposition}[\cite{P-3}]\label{Proposition 2.2 Trieste}
The space of Fox $k$-colorings is preserved by $k$-moves.
\end{proposition}
\begin{proof}  Figure 2.3 illustrates the bijection between $col_k(D)$ and $col(m_k(D))$ where 
$m_k(D)$ is obtained from $D$ by a $k$-move. This bijection is an isomorphism of groups 
$Col_k(D)$ and $Col_k(m_k(D))$
\end{proof}

\ \\
\centerline{\psfig{figure=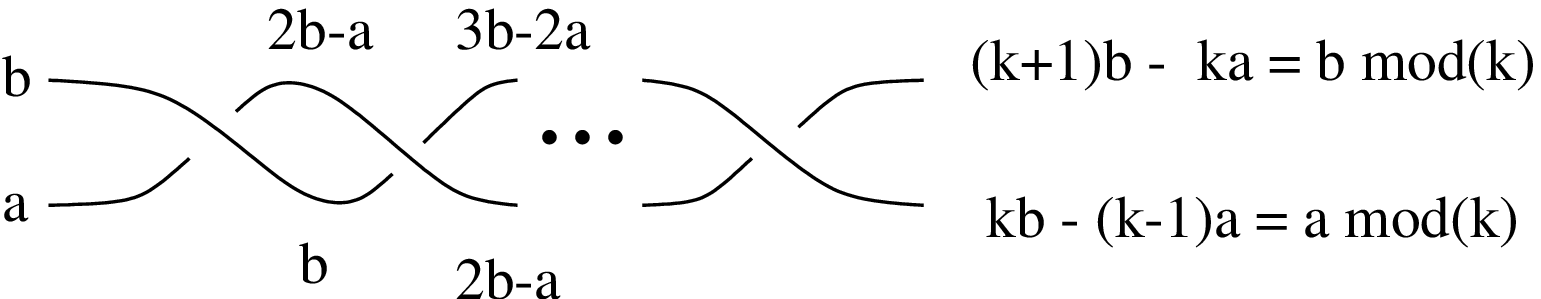,height=2.6cm}}
\centerline{Fig. 2.3; from $(b,a)$ to $(k(b-a)+b, k(b-a)+a)$}
\ \\

The following properties of k-colorings, are a straightforward generalizations from 
3-colorings and can be proved in a similar way. However, an elementary proof
of the part (c) is, as before, more involved and the simplest proof (not involving double branched covers), 
I am aware of, requires an interpretation
of k-colorings using the Goeritz matrix \cite{Goe,Gor-1,P-5} or use of Lagrangian tangles (see below).
\begin{lemma}\label{Lemma 2.3 Trieste}
\begin{enumerate}
\item
[(a)] $col_k(L)$ is a divisor of a power of $k$ and for a
link with $b$ bridges, $col_k(D)$ divides $k^b$. More precisely.
 $Col_k(L)$ is a subgroup of $Z_k^b$.
\item
[(b)] $col_k(L_1) col_k(L_2)=k(col_k(L_1\# L_2))$ (notice that our yesterday's proof works only 
for odd $k$ as we use the fact that $2$ is invertible in $Z_k$), 
\item
[(c)] Consider $k+1$ diagrams $L_0, L_1, ..., L_{k-1},L_{\infty}$; see Fig.~2.4 .
If $k$ is a prime number then among the $k+1$ numbers
$col_k(L_0), col_k(L_1),...,col_k(L_{k-1})$ and $col_k(L_{\infty})$
$k$ are equal one to another and the $(k+1)$'th is $k$ times bigger.
\end{enumerate}
\end{lemma}
\ \\

\centerline{\psfig{figure=L0L1L2L3Linf.eps,height=1.5cm}}
\begin{center}
Fig. 2.4
\end{center}

Notice, that (c) can be interpreted as follows: \\
Let $k$ be a prime number and $col'_k(L)=(-1)^{col_k (L)}col_k (L)$ then 
$$col'_k(L_0)+ col'_k(L_1)+...+ col'_k(L_{k-1}) + col'_k(L_{\infty}) =0, $$ 
a skein relation of $k+1$ terms often called $(k,\infty)$ skein relation.

\begin{example}\label{Example 2.3}
\begin{enumerate}
\item
[(i)] For the figure eight knot, $4_1$, one has $col_5(4_1)=25$, so the figure
eight knot is a nontrivial knot; compare Figure 2.5.
\item
[(ii)] For the knot $5_2$ we have $col_7(5_2)=49$ (more precisely $Col_7(5_2)=\Z_7^2$);  
compare Figure 2.5.
\end{enumerate}
\end{example}

\centerline{\psfig{figure=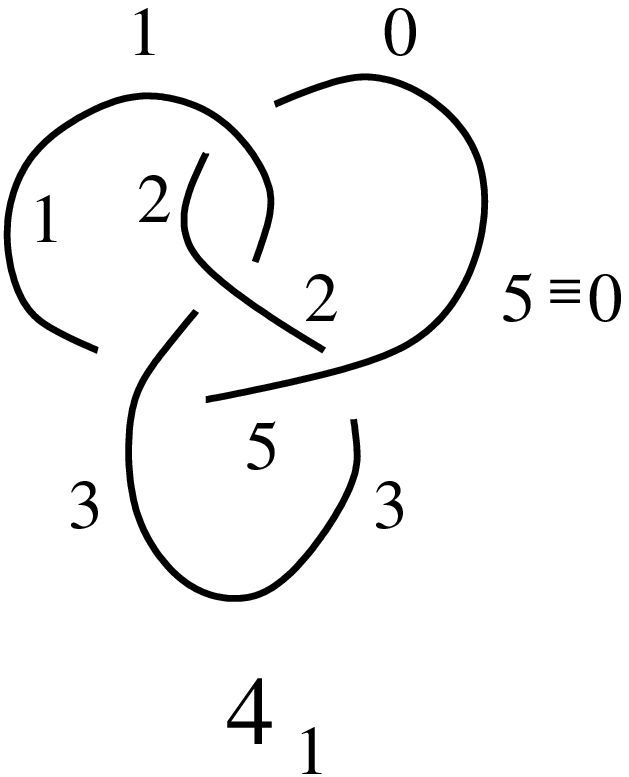,height=5.9cm} \ \ \ \ \ \ \ \psfig{figure=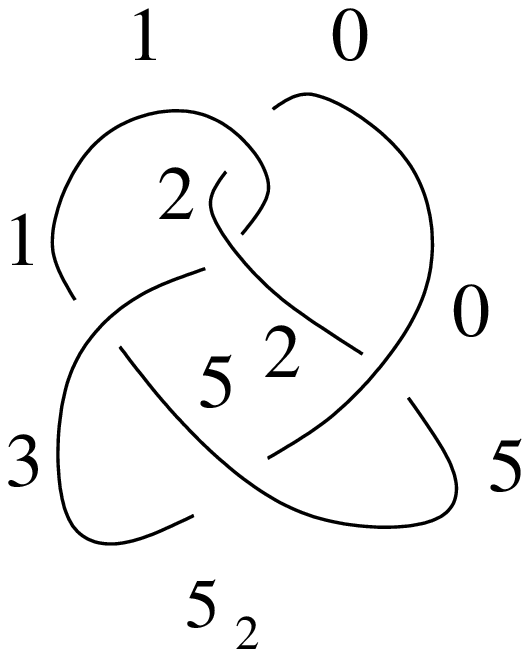,height=5.9cm}}
\begin{center}
Fig. 2.5
\end{center}

Let us look closer at the observation that a $k$-move
preserves the space of Fox $k$-colorings.
One should consider a general  {\it rational moves}, that is,
a rational $\frac{p}{q}$-tangle of Conway is substituted in
place of the identity tangle\footnote{The move was first
considered by J.M.Montesinos\cite{Mo-2}; compare also Y.Uchida \cite{Uch}.}.
The important observation for us is that
$Col_{p}(D)$ is preserved by $\frac{p}{q}$-moves. Fig.2.6 
illustrates the fact that $Col_{13}(D)$ is unchanged by
a $\frac{13}{5}$-move.
\\
\ \\
\centerline{\psfig{figure=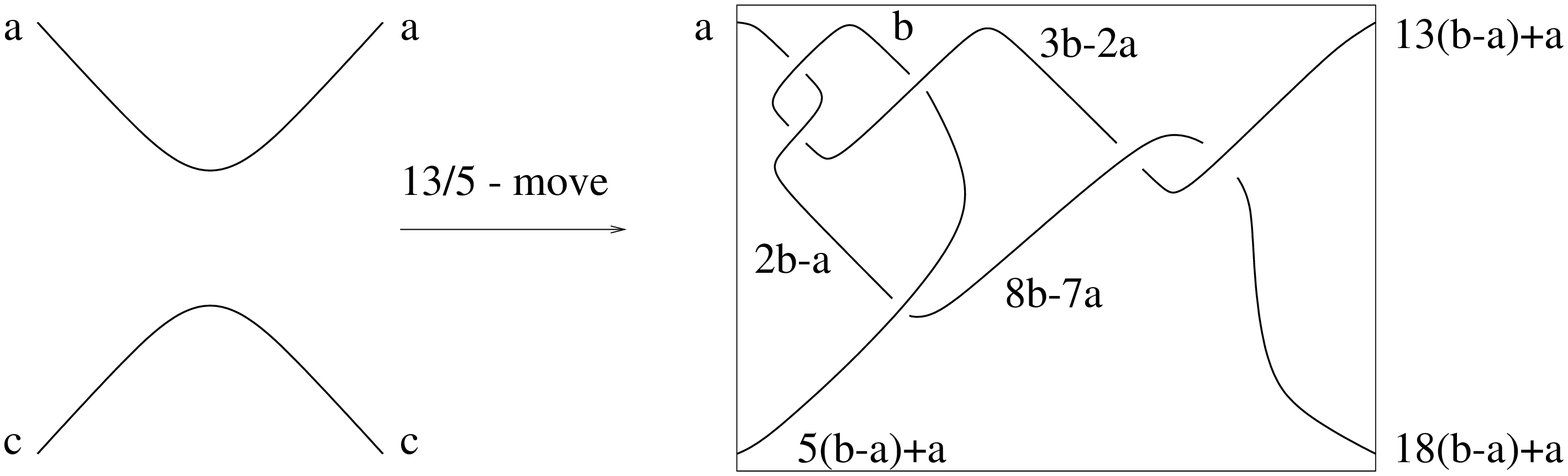,height=4.3cm}}
\ \\ \ \\

\centerline{\psfig{figure=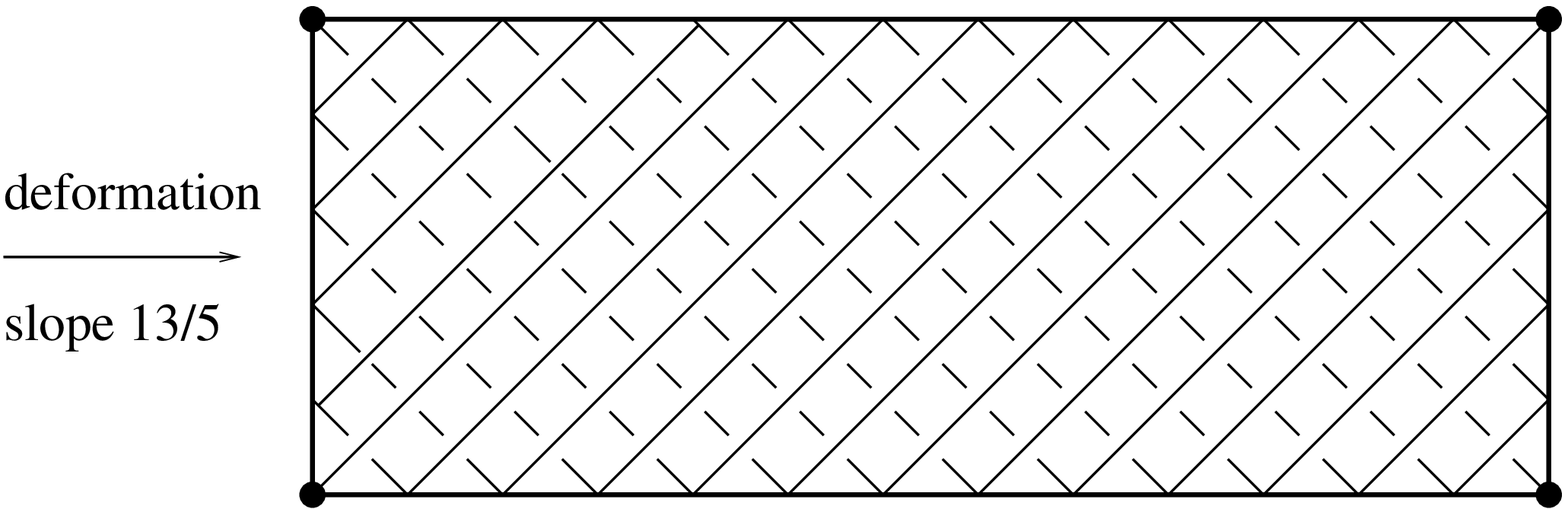,height=3.6cm}}\ \\

\centerline{Fig. 2.6; $\frac{13}{5}$-move and $\frac{13}{5}$ tangle in 
Conway's and pillow case forms}
\ \\ \ \\ 
We just have heard about the Conway's  classification
of rational tangles at the Lou's and Sofia's talks, so  I only  
briefly sketch definitions and notation.
The 2-tangles shown in Figure 2.7 are called rational tangles
with Conway's notation $T(a_1,a_2,...,a_n)$.
A rational tangle is the
$\frac{p}{q}$-tangle if $\frac{p}{q} =
a_n + \frac{1}{a_{n-1}+...+\frac{1}{a_1}}$.\footnote{$\frac{p}{q}$
is called the slope of the tangle and can be easily
identified with the slope of the meridian disk of the solid torus
being the branched double cover of the rational tangle.} Conway
proved that two rational tangles are ambient isotopic
(with boundary fixed) if and only if their slopes
are equal (compare \cite{Kaw}).
\\
\ \\
\centerline{\psfig{figure=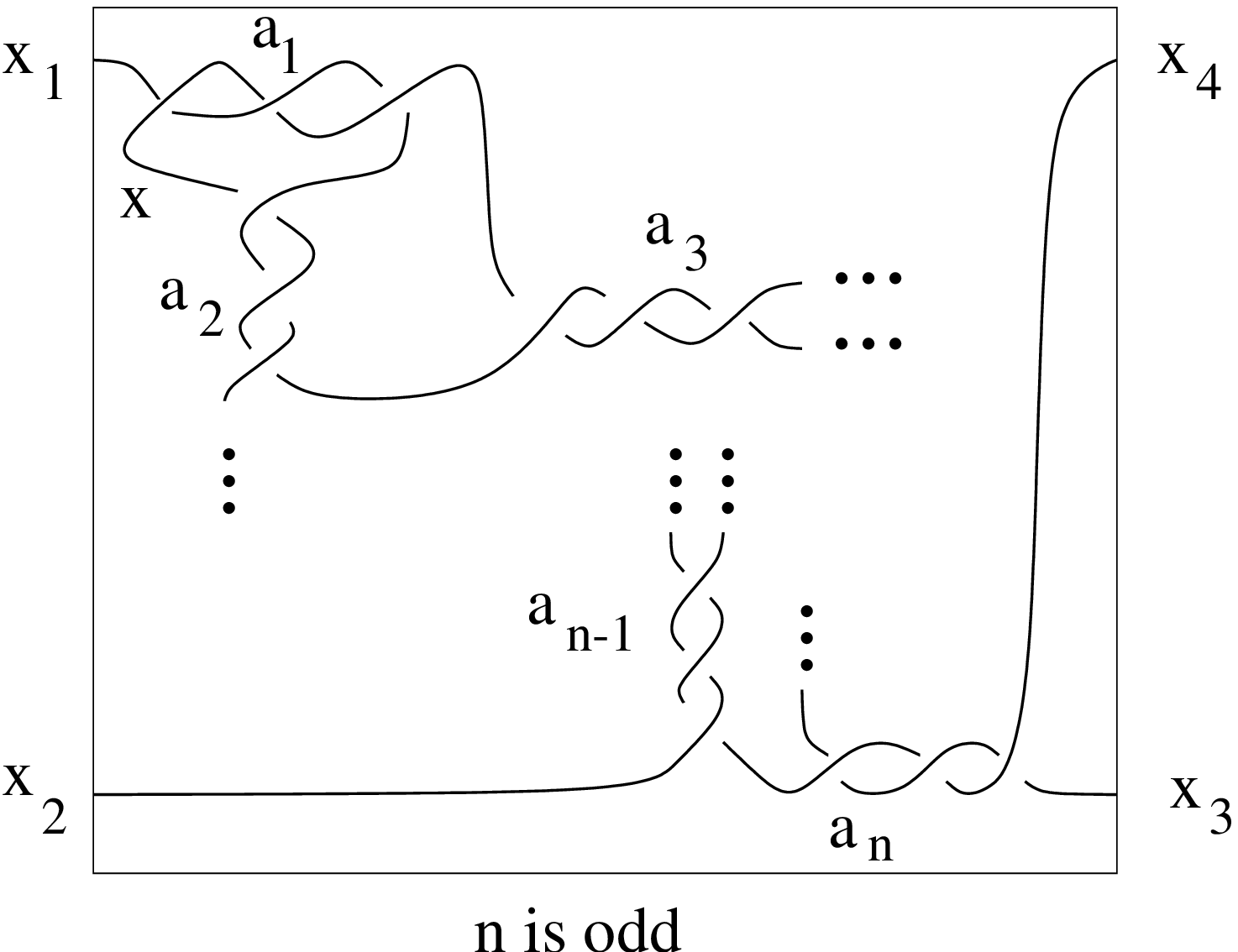,height=4.9cm}\ \ \
\psfig{figure=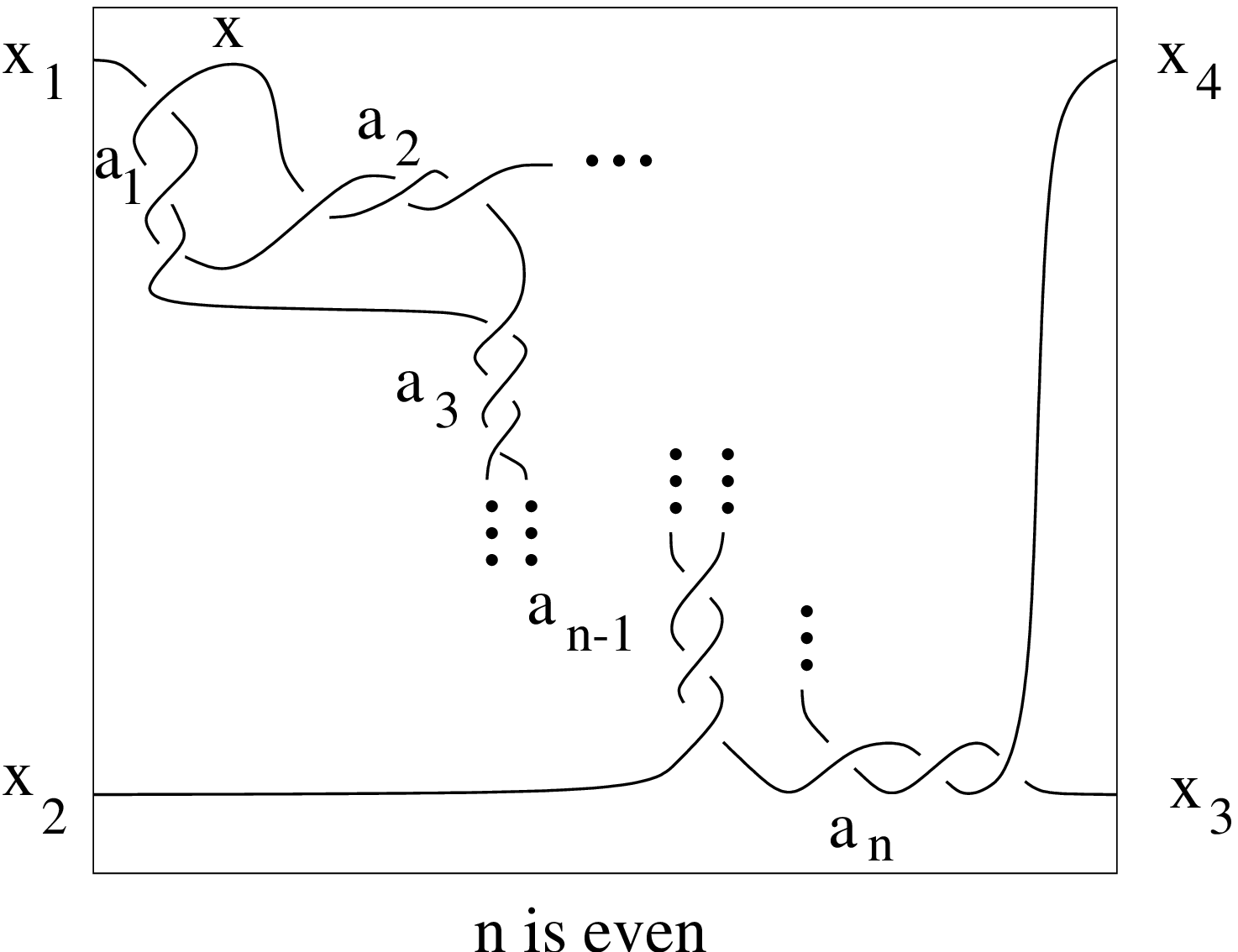,height=4.9cm}}
\centerline{Fig. 2.7}

For a Fox coloring of a rational $\frac{p}{q}$-tangle with
boundary colors $x_1,x_2,x_3,x_4$ (Fig.2.5), one has
$x_4-x_1 = p(x-x_1)$, $x_2-x_1=q(x-x_1)$ and $x_3 =
x_2 + x_4 -x_1$. If a coloring is nontrivial ($x_1\neq x$)
then $\frac{x_4-x_1}{x_2-x_1} = \frac{p}{q}$ as
has been explained in the talk by Lou Kauffman.

\begin{corollary}\label{Corollary 2.5} 
$\frac{p}{q}$-move on a link or a tangle is preserving the group of $p$-colorings.
\end{corollary}

\subsection{Symplectic structure on Fox Colorings,  Lagrangian tangles}

The usefulness of the symplectic structure in the knot theory,
was probably first observed by R.~Fox in his review of the A.~Plans
paper of 1953 \cite{Pla}. 
In this part we follow \cite{DJP}  showing how to define a symplectic form on the space of
Fox colorings of the boundary of
n-tangles so that every tangle corresponds to a Lagrangian (in
the case of a field of colors) or a virtual Lagrangian (for PID) of
the symplectic structure (that is, a subspace of a maximal
dimension on which the form vanishes). 
Inversely, for a field $R=Z_p$, $p>2$, every Lagrangian can
be realized by a tangle. It does not hold for $Z_2$ and $n>3$.\footnote{In \cite{DJP} we 
draws from the construction several far fetching conclusions: first, it allows us 
to understand the space of colored tangles as a Tits building. Second, it provides 
 applications to 3-manifold topology.
In particular, we show that our symplectic space is related (via double
branched cover) to the symplectic structure on homology on a surface
(with the symplectic form given by the intersection number). It relates
our results with a known fact that 3-manifolds yield Lagrangians
in $H_1(\partial M;Q)$.  One application is to use Lagrangians
to find obstructions for embedding n-tangles into links.
Rotation of a tangle yields an isometry of our symplectic space, and we
analyze invariant subspaces of the map, in particular we look for
invariant Lagrangians of the rotation by $2\pi/n$ (along $z$-axis).
We use our analysis to answer, partially the question whether
rotation of a link (as described in \cite{APR}) preserves the homology
of the double branch cover of $S^3$ with the link as branching set.}

\subsection{Alternating form on colorings of a tangle boundary}

We work with modules over a commutative ring with identity, $R$.
We concentrate our attention on the finite field $R=Z_p$.

Consider $2n$ points on a circle (or a square, Fig. 2.8). Let the ring $R$
be treated as a set of colors (e.g. a field ${\bf Z}_p$).

For $R={\bf Z}_p$ the colorings of
$2n$ points form a linear space $V={\bf Z}_{p}^{2n}$. Let $e_1,\ldots
, e_{2n}$ be its basis, $e_i=(0,\ldots , 1,\ldots ,0)$, where 1
occurs in 
\ \\
\begin{center}
\centerline{\psfig{figure=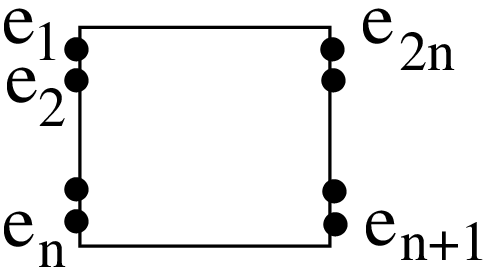}}
\end{center}
\begin{center} Fig. 2.8 \end{center}
\ \\
 the $i$-th position. Let $V'={\bf Z}_{p}^{2n-1}\subset{\bf
Z}_p^{2n}$ be the subspace of vectors $\sum a_i e_i$ satisfying $\sum (-1)^i
a_i=0$ (the alternating condition).
Consider the basis $f_1,\ldots , f_{2n-1}$ of ${\bf
Z}_{p}^{2n-1}$ where $f_k=e_k+e_{k+1}$. We can also introduce the vector $f_{2n}=e_{2n} + e_{1}$ and then 
$f_{2n}= f_1-f_2+f_3\pm... -f_{2n-2}+ f_{2n-1}$.
 Consider an alternating form\footnote{That is for any $a\in V$ one has $\phi(a,a)=0$. From this 
anti-symmetry follows ($\phi(a,b)= -\phi(b,a)$).}
 $\phi$ on ${\bf Z}_{p}^{2n-1}$ of nullity 1 given by the matrix

$$\phi = \left( \begin{array}{cccccccc} 
0 & 1 & 0 & 0 & \ldots & 0 & 0 & 0 \\
 -1 & 0 & 1 & 0 &\ldots & 0 & 0 & 0 \\ 
0 & -1 & 0 & 1 &\ldots & 0 & 0 & 0\\
\ldots & \ldots & \ldots & \ldots &\ldots & \ldots & \ldots & \ldots \\ 
\ldots & \ldots & \ldots & \ldots & \ldots &  \ldots & \ldots & \ldots \\
0 & 0 & 0 &  0 & \ldots & -1 & 0 & 1\\
 0 & 0 & 0 &  0 & \ldots & 0& -1 & 0
\end{array} \right)$$

that is,
 \vspace{1mm}
   \renewcommand{\arraystretch}{3}
   $$\phi (f_i, f_j) =\left\{
   \begin{array}{lr}
   0 &
   {\rm if}\ |j-i|\neq 1 \\
  1 &  {\rm if\ }\ j=i+1\\ -1 &
   {\rm if}\ j=i-1. \\
   \end{array}
  \right .$$
   \par
   \vspace{2mm}

An alternating form of nullity one is called a pre-symplectic form. 
A pre-symplectic form on ${\bf Z}_{p}^{2n-1}$ leads to a symplectic (i.e. alternating, non-degenerated) 
form on ${\bf Z}_{p}^{2n-2}$ as follows:\\
The vector $e_1+ e_2+ \ldots + e_{2n}$
($=f_1+f_3+\ldots + f_{2n-1}=f_2+f_4+\ldots + f_{2n})$
is $\phi$-orthogonal to any other vector.
If we consider ${\bf Z}_{p}^{2n-2}={\bf Z}_{p}^{2n-1}/{\bf
Z}^{tr}_{p}$, where the subspace ${\bf Z}^{tr}_{p}$ is generated by
$e_1+\ldots + e_{2n}$, that is, ${\bf Z}^{tr}_{p}$ consists
of monochromatic (i.e. trivial) colorings,
then $\phi$ descends to a symplectic form $\hat\phi$ on ${\bf
Z}_{p}^{2n-2}$. Now we can analyze the isotropic subspaces of
$({\bf Z}_{p}^{2n-2},\hat\phi)$, that is, subspaces on
which  $\hat\phi$ is $0$
($W\subset {\bf Z}_{p}^{2n-2}, \phi (w_1,w_2)=0$ for $w_1,w_2\in W$).
The maximal isotropic ($(n-1)$-dimensional)
subspaces of ${\bf Z}_{p}^{2n-2}$ are
called Lagrangian subspaces (or maximal totally degenerated subspaces)
and there are $\prod_{i=1}^{n-1}(p^i+1)$ of them. We use the term 
pre-Lagrangian for a maximal totally degenerated subspace of ${\bf Z}_{p}^{2n-1}$. 
Of course, ${\bf Z}^{tr}_{p}$ lies in every pre-Lagrangian. Lagrangians in ${\bf Z}_{p}^{2n-2}$ 
are $(n-1)$-dimensional and pre-Lagrangians in ${\bf Z}_{p}^{2n-1}$
are $n$-dimensional \cite{O'M}. 

Let $\psi=\psi_T :Col_p (T)\rightarrow{\bf Z}_{p}^{2n}$ be the homomorphism which sends colorings of $T$ into 
colorings of boundary points of the tangle. Our local
condition on Fox colorings (Fig.2.1) guarantees that for any n-tangle $T$,
$\psi (Col_p (T))\subset {\bf Z}_{p}^{2n-1}$.\footnote{We checked it before for a ring $R$ in which
$2$ is not a zero divisor; the general case follows from considerations given later.}
 Furthermore,
the space of trivial colorings, ${\bf Z}^{tr}_{p}$, always lies in $Col_p (T)$. 
The quotient space $Col_p (T)/{\bf Z}^{tr}_{p}$ is called the reduced space of Fox colorings 
and denoted by $Col^{rd}_p (T)$.
Thus $\psi$ descends to $\hat\psi :Col^{rd}_p (T) \rightarrow
{\bf Z}_{p}^{2n-2}={\bf Z}_{p}^{2n-1}/{\bf Z}_{p}$. Now we have
a fundamental question:
which subspaces of
${\bf Z}_{p}^{2n-2}$ are yielded by $n$-tangles?
We answer this question below.
\begin{theorem}\label{Theorem 2.6 Trieste}
$\hat\psi (Col^{rd}_p (T))$ is a Lagrangian
subspace of ${\bf Z}_{p}^{2n-2}$ with the symplectic form $\hat\phi$. 
In particular, $dim (\hat\psi (Col^{rd}_p (T)))=n-1$.\ 
Equivalently, $\psi (Col_p (T))$ is a pre-Lagrangian
subspace of ${\bf Z}_{p}^{2n-1}$ with the alternating form $\phi$. 
In particular, $dim (\psi (Col_p (T)))=n$.
\end{theorem}
A natural question would be whether every Lagrangian subspace
can be realized by a tangle. The answer is negative for
$p=2$, but for $p>2$ we have
\begin{theorem}\label{Theorem 2.7 Trieste}
For $p$ an odd prime number every Lagrangian subspace of ${\bf Z}_{p}^{2n-2}$ can be realized 
by a tangle, in fact, by an $n$-rational tangle\footnote{An $n$-rational tangle is an n-tangle having 
 presentation so it is often called an $n$-bridge tangle.}.
\end{theorem}
Theorem \ref{Theorem 2.7 Trieste} follows from the work of J.~Assion \cite{As} (for $p=3$),
and B.~Wajnryb \cite{Wa-1,Wa-2} (for $p>2$). Wajnryb constructs the natural epimorphism 
from the odd braid group $B_{2n+1}$ to the symplectic group ${\rm Sp}(n,p)$, that is, the group of isometries 
of the symplectic space ${\bf Z}_{p}^{2n}$.

As a corollary to Theorems \ref{Theorem 2.6 Trieste} and \ref{Theorem 2.7 Trieste} we obtain 
a fact which was considered difficult before, even for 2-tangles.
\begin{corollary}\label{Corollary 2.8 Trieste}
For any $p$-coloring of a tangle boundary satisfying
the alternating property (i.e., is an element
of ${\bf Z}_{p}^{2n-1}$) there is an $n$-tangle and
its $p$-coloring yielding the given coloring on the boundary.
In other words: ${\bf Z}_{p}^{2n-1} = \bigcup_T \psi_T(Col_p(T))$.
Furthermore, the space $\psi_T(Col_p(T))$ is $n$-dimensional.
\end{corollary}

In \cite{DJP} we give short, high-tech  proof of Theorems 2.6 and 2.7. Here we provide a longer but 
elementary proof based on the presentation of an n-tangle as a 
tangle with $N$ maxima and $N-n$ minima as in Figure 1.13.
The proof of Theorem \ref{Theorem 2.6 Trieste} is straightforward: we check 
that the theorem holds for a trivial n-tangle, that it holds when we add crossings,
and finally that it is still valid after applying minima. On the way we show that 
Theorem \ref{Theorem 2.7 Trieste} follows from the second part of the proof (the slogan will be that 
braid-like transvections generate a symplectic group over ${\bf Z}_p$, $p>2$; as follows from Wajnryb \cite{Wa-2}).\\
{\bf Step 0:} Consider the trivial tangle, $T_0$,
 in which the point $v_{2i-1}$ is connected to $v_{2i}$; see Figure 2.9. Clearly $\psi(Col_p(T_0))$ 
is the n-dimensional subspace of ${\bf Z}_{p}^{2n}$ generated by $e_{2i-1}+e_{2i}$ ($1 \leq i \leq n)$, and 
thus it is a pre-Lagrangian in ${\bf Z}_{p}^{2n-1}$ generated by $f_1,f_3,...,f_{2n-1}$. In effect,
$\hat\psi(Col^{rd}_p(T_0))$ is the Lagrangian subspace of ${\bf Z}_{p}^{2n-2}$ generated by 
$f_1,f_3,...,f_{2n-3}$.

\ \\
\begin{center}
\centerline{\psfig{figure=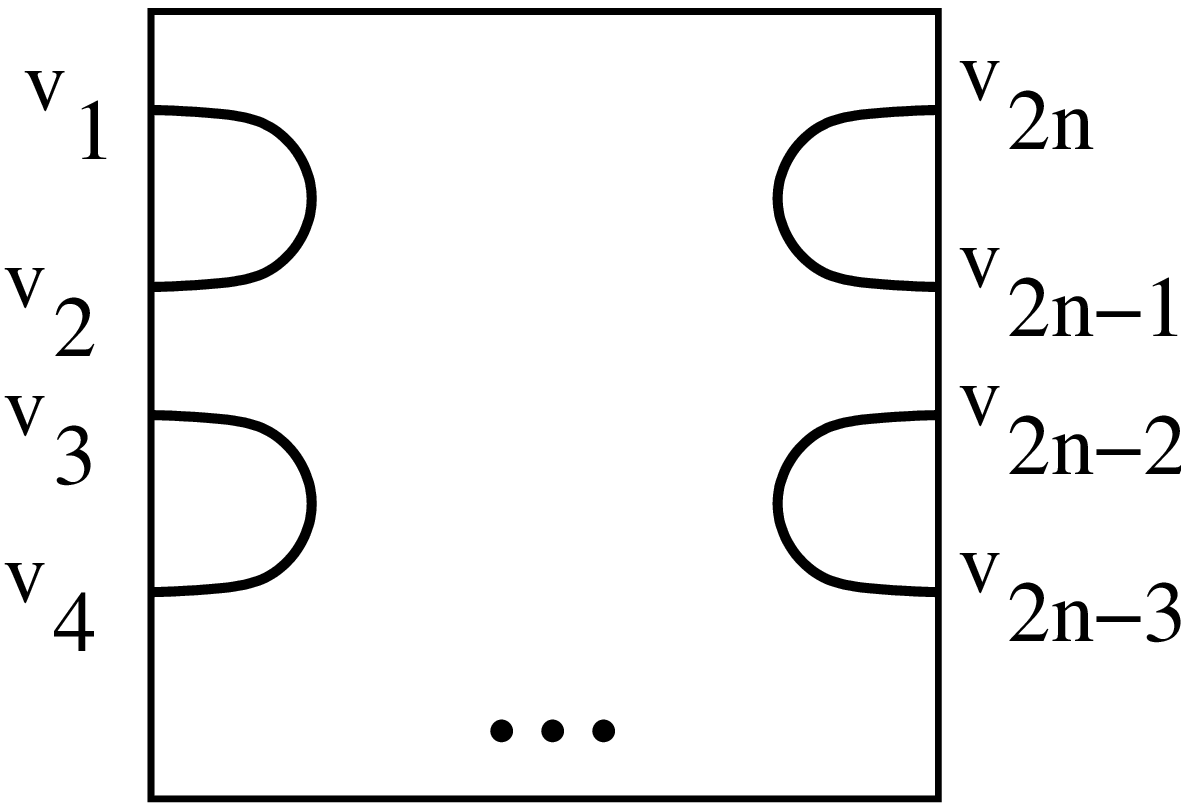,height=3.1cm}}
\end{center}
\begin{center} Fig. 2.9 Trivial $n$-tangle, $T_0$ \end{center}
\ \\

{\bf Step 1:} We show here that if $\hat\psi_T(Col^{rd}_p(T))$ is a Lagrangian in ${\bf Z}_{p}^{2n-2}$ 
then $\hat\psi(Col^{rd}_p(T'))$ is also a Lagrangian where $T'$ is a tangle obtained from $T$ by adding
one  crossing to it (without loss of generality we can assume that the crossing is between arcs from $v_{2n}$ 
and $v_{2n-1}$, see Figure 2.10; the relevant observation is that the rotation of the tangle by $\frac{2\pi}{2n}$ is an 
isometry, that is $\phi(f_i,f_j)=\phi(f_{i+1},f_{j+1})$, where indices are taken modulo $2n$).  
\ \\
\begin{center}
\centerline{\psfig{figure=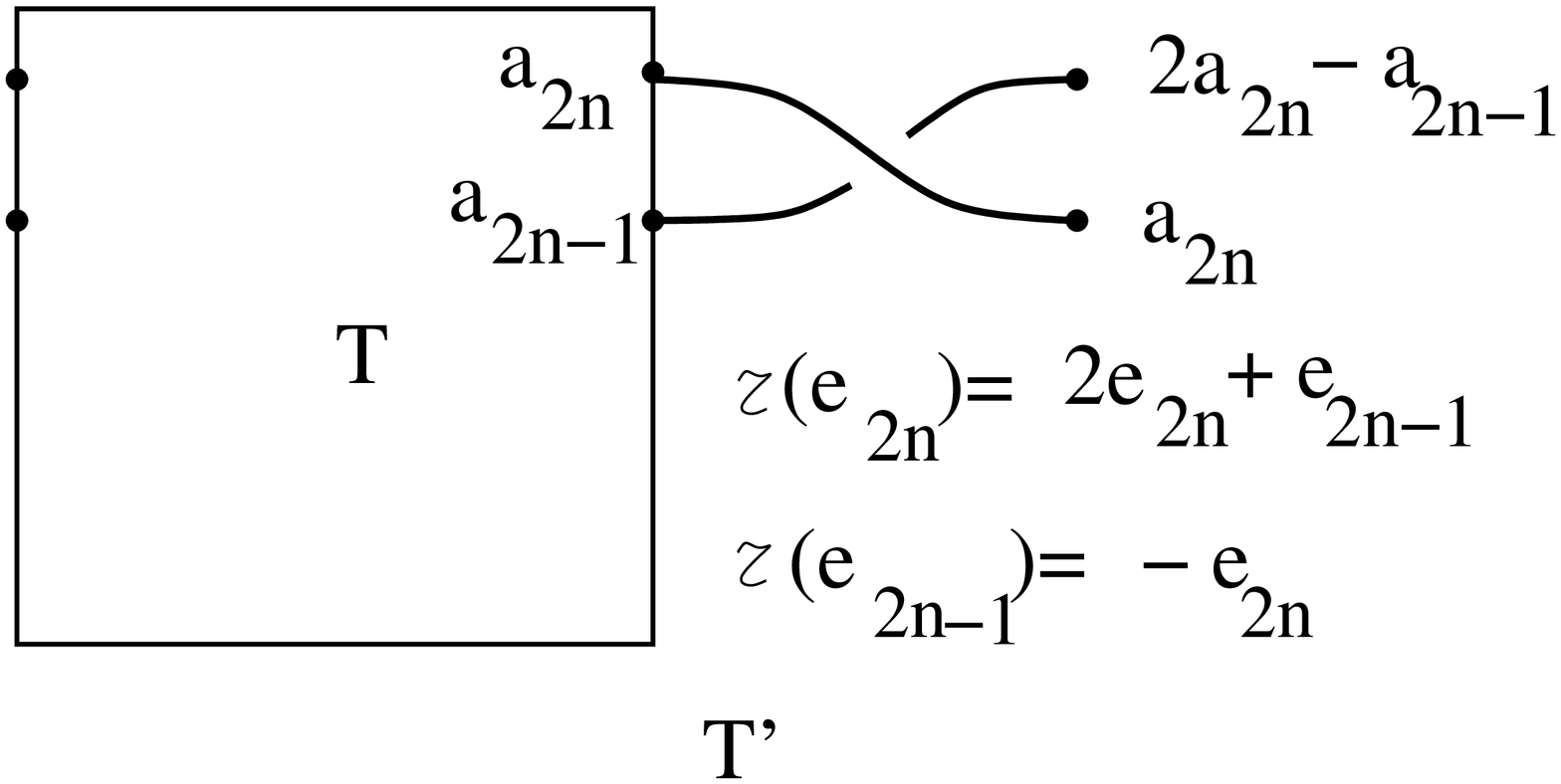,height=4.1cm}}
\end{center}
\begin{center} Fig. 2.10; adding a crossing to $T$ induces isometry on $({\bf Z}_{p}^{2n-1},\phi)$ \end{center}
\ \\
Moving from a coloring of the boundary of $T$ to the boundary of $T'$ induces the linear map
$\tau: {\bf Z}_{p}^{2n} \to {\bf Z}_{p}^{2n}$ given by:\\ 
$\tau(e_{2n})= 2e_{2n} + e_{2n-1}$, $\tau(e_{2n-1})= -e_{2n}$ and $\tau(e_{i}))= e_i$ otherwise. 
$\tau$ sends an alternating sum to an alternating sum so it preserves the subspace ${\bf Z}_{p}^{2n-1}$ 
on which the alternating form $\phi$ is defined. In the basis $f_1,f_2,...,f_{2n-1}$ it is defined by:
$$\tau(f_{2n-2}) = \tau(e_{2n-2}+e_{2n-1})= e_{2n-2}-e_{2n}= f_{2n-2}-f_{2n-1}$$
$$\tau(f_{2n-1}) = \tau(e_{2n-1}+e_{2n})= e_{2n} + e_{2n-1} = f_{2n-1}$$
$$\tau(f_{2n}) = \tau(e_{2n}+e_{1})= 2e_{2n} +  e_{2n-1} +e_{1}= f_{2n}+f_{2n-1}$$  
$$ \tau(f_{i}) = \tau(f_{i})\ \ otherwise.$$ 
In summary, we get $ \tau(f_{i}) = f_i - \phi(f_i,f_{2n-1})f_{2n-1}$. Such a linear map 
is called a symplectic transvection with respect to vector $f_{2n-1}$ and denoted by $\tau_{f_{2n-1}}$.
Generally the transvection $\tau_{b}(a)= a - \phi(a,b)b$ is an isometry with respect to the form $\phi$
(for completeness here is the check: 
$$\phi(\tau_{b}(a_1), \tau_{b}(a_2))= \phi(a_1-\phi(a_1,b)b, a_1-\phi(a_2,b)b)= $$
 $$ \phi(a_1,a_2) - \phi(a_1, \phi(a_2,b)b) - \phi(\phi(a_1,b)b,a_2) =  \phi(a_1,a_2).)$$
Notice that if we change the crossing in Figure 2.10 to its mirror image then $\tau$ is replaced 
by $\tau^{-1}$ with $\tau^{-1}(f_i)= f_i + \phi(f_i,f_{2n-1})f_{2n-1}$. 

For us it is important that transvection, as an isometry, is sending pre-Lagrangians to pre-Lagrangians 
and Lagrangians to Lagrangians.

{\bf Step 2:}\footnote{It is related to the ``contraction lemma" described by Turaev \cite{Tur}, p.180.}
 Consider a minimum (here right cup, see Figure 2.11).

\ \\
\begin{center}
\centerline{\psfig{figure=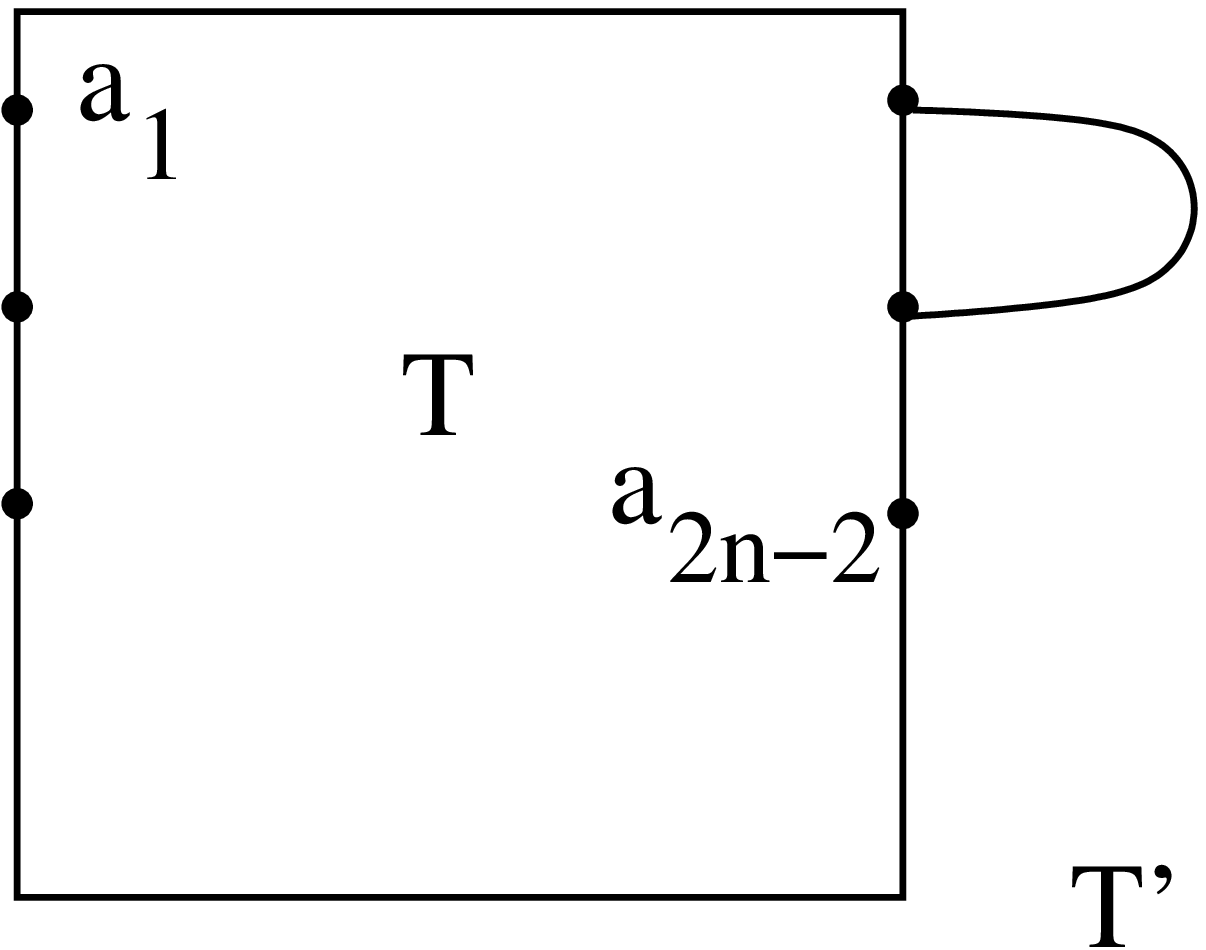,height=4.1cm}}
\end{center}
\begin{center} Fig. 2.11; adding a minimum (right cup)  \end{center}
\ \\

Initially let us consider adding a right cup in general, without assuming that 
$\psi_T(Col_p(T))$ is a pre-Lagrangian in ${\bf Z}_{p}^{2n-1}$, (this may be useful in 
a less restricted setting of virtual or welded tangles).\\
Consider the linear space $V$ with a basis $\{e_1,...,e_{2n}\}$ (corresponding to 
${\bf Z}_{p}^{2n}$ and ${\bf Z}_p$ colorings of the boundary of a $n$-tangle) and consider the right 
cup of Figure 2.11. We analyze induced ${\bf Z}_p$-colorings of the boundary of a $(n-1)$-tangle in two steps:\\
(I) Let $F$ be a subspace of $V$:
\begin{enumerate}
\item[(1)] We consider the subspace $F_1$ of $F$ defined by $F_1=\{ a=\sum_{i=1}^{2n}a_i \in F\ | \ a_{2n-1}=a_{2n}\}$. 
We have two cases for the dimension of $F_1$.
\begin{enumerate}
\item[(i)] $dim (F_1) = dim (F)$. This is the case iff $F_1\subset span\{e_1,...,e_{2n-2},f_{2n-1}\}$, where 
$f_{2n-1}=e_{2n-1}+ e_{2n}$.
\item[(ii)] $dim (F_1) = dim (F) -1$.  This is the case iff there is $a\in F$ such that $a_{2n-1}\neq a_{2n}$.
\end{enumerate}
\item[(2)] We consider the projection $p:V \to W=span\{e_1,...,e_{2n-2}\}$ 
(here $p(e_{2n-1})=p(e_{2n}) = 0$ and $p(e_i)=e_i$ for $1 \leq i \leq 2n-2$.). Let $F_2=p(F_1)\subset W$.
We have two cases for the dimension of $F_2$:
\begin{enumerate}
\item[(i)] $dim(F_2)=dim(F_1)$. This holds iff $f_{2n-1}$ is not in $F_1$.\\
\item[(ii)] $dim(F_2)=dim(F_1)-1$ iff $f_{2n-1} \in F_1$.
\end{enumerate}
\end{enumerate}
We show that both (1)(i) and (2)(i) cannot hold if $F$ is a pre-Lagrangian. Similarly 
(1)(ii) and 2(ii) cannot hold for such an $F$. 

If elements of $F$ satisfy the alternating condition ($a\in F \Rightarrow \sum_{i=1}^{2n}(-1)^ia_i=0$),
then (1) can be reformulated as: 
\begin{enumerate}
\item[(1)(i')] $dim (F_1) = dim (F)$ iff $F_1\subset span\{f_1,...,f_{2n-3},f_{2n-1}\}$, \\
\item[(1)(ii')] $dim (F_1) = dim (F)-1$ iff there is $a\in F$ such that 
$a =f_{2n-2}+ v$, $v\in span\{f_1,...,f_{2n-3},f_{2n-1}\}$.
\end{enumerate}
Let $V'$ be the subspace of $V$ of elements satisfying the alternating condition, thus $V'$ is generated by 
$f_1,...,f_{2n-1}$. Assume that  $F$ is a pre-Lagrangian in $V'$. Then 1(i') and 2(i) cannot 
hold as in that case  $F$ is not a pre-Lagrangian -- it is not maximal: adding $f_{2n-1}$ still gives
a totally degenerated space. .\\ 
Similarly, if 1(ii') and 2(ii) hold then
$a\in F$ and $f_{2n-1}\in F$ but \\
$\phi(a,f_{2n-1})= \phi(f_{2n-2},f_{2n-1})=1$, so $F$  could not be a totally degenerate space.\\

Thus we proved that $F_2$ is $n-1$ dimensional in $W$. It is also a totally degenerated space 
(as an embedding of $W'$ in $V'$ is an isometry: $W'$ is a subspace of $W$ satisfying the
alternating condition, so it has a basis $f_1$,...$f_{2n-2}$). Thus $F_2$ is a pre-Lagrangian in $W'$.
The proof of Theorem \ref{Theorem 2.6 Trieste} is completed.

As we mentioned before, Theorem \ref{Theorem 2.7 Trieste} follows from the result of Wajnryb that 
the symplectic group is generated by braid-like generators in which braid generators act as transvections \cite{Wa-2}. 
I our situation it means that adding crossings to a tangle allows us to realize any symplectic map on 
the symplectic space ${\bf Z}_p^{2n-2}$ of boundary coloring. In particular any Lagrangian is an 
image of the Lagrangian $span\{f_1,f_3,...,f_{2n-2}\}$ associated to the trivial n-tangle $T_0$. 
\ \\
\ \\
\centerline{\psfig{figure=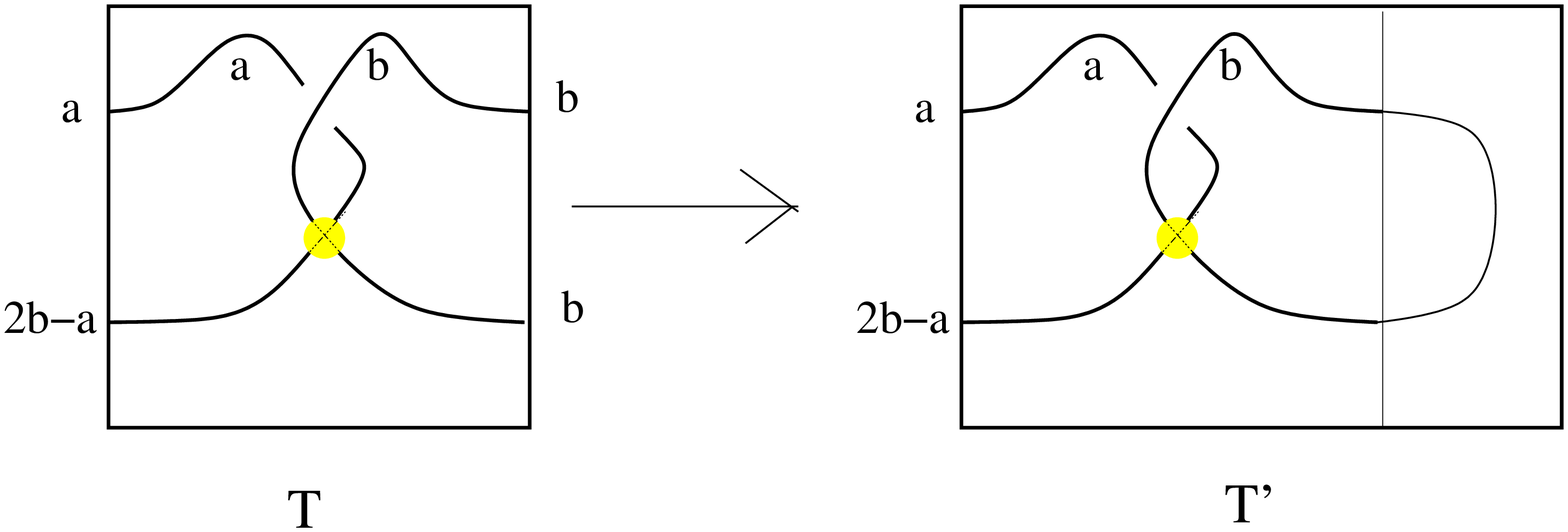,height=2.6cm}}
\centerline{Fig. 2.12; non-classical colorings of virtual tangles}
\ \\

On Figure 2.12 we have an example of a virtual 1-tangle $T'$ such that $\psi(Col_p(T'))={\bf Z}_p^{2}$.
Combining $n$ tangles $T'$ together we get a virtual $n$-tangle $T^{(n)}$ with $n$ virtual crossings 
and $\psi(Col_p(T^{(n)}= {\bf Z}_p^{2n}$. Combining $T'$ tangles and trivial tangles we can 
get a virtual tangle $T$ with $dim (\psi(Col_p(T)))$ any number between $n$ and $2n$. Can we get 
dimension smaller from $n$? In particular, is there a virtual 2-tangle such that boundary coloring is always 
monochromatic?

 Let me complete this presentation by mentioning two generalizations
of the Fox $k$-colorings.

In the first generalization we consider any commutative
ring with the identity in place of ${\bf Z}_k$. We construct $Col_RT$
in the same way as before with the relation at
each crossing, Fig.2.1,
having the form $c=2a-b$ in $R$.
The skew-symmetric form $\phi$ on $R^{2n-1}$, the symplectic
form $\hat\phi$ on $R^{2n-2}$ and the homomorphisms $\psi$ and
$\hat\psi$ are defined in the same manner as before.
Theorem  2.4 generalizes as follows (\cite{DJP}):
\begin{theorem}\label{2.7}
Let $R$ be a Principal Ideal Domain (PID) then,
$\hat\psi (Col_R T/R)$ is a virtual Lagrangian
submodule of $R^{2n-2}$ with the symplectic form $\hat\phi$.
That is $\hat\psi (Col_R T/R)$ is a finite index submodule
of a Lagrangian in $R^{2n-2}$.
\end{theorem}

The second generalization leads to racks and quandles \cite{Joy,F-R}
but we restrict our setting to the abelian case -- Alexander-Burau-Fox
colorings\footnote{The related approach was first outlined in the letter
of J.W.Alexander to O.Veblen, 1919 \cite{A-V}. Alexander was probably
influenced by P.Heegaard dissertation, 1898, which he reviewed for
the French translation \cite{Heeg}. Burau was considering a braid
representation but locally his relation was the same as that of Fox.
According to J.Birman, Burau learned of the representation
from Reidemeister or Artin \cite{Ep}, p.330.}.
An ABF-coloring uses colors from a ring, $R$, with an invertible
element $t$ (e.g. $R= \Z[t^{\pm 1}]$). The relation in Fig.2.1
is modified to the relation $c=(1-t)a + tb$ in $R$ at each crossing
of an oriented link diagram; see Fig. 2.13.
\ \\
\ \\
\centerline{\psfig{figure=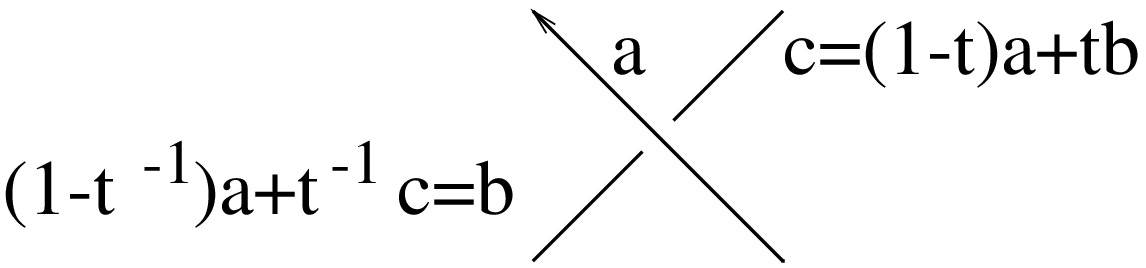,height=2.6cm}}
\centerline{Fig. 2.13}

The space $R^{2n-2}$ has a natural Hermitian structure \cite{Sq},
 one can also find a symplectic structure and one can prove Theorem 2.7
in this setting \cite{DJP}.

\section{Conclusion}\label{Section 3}
I hope that our snapshot of knot theory will inspire you to consider ideas described 
in the last two days. I am sure you are already asking: what about other $n$-move conjectures?
Why should we use only abelian groups? Can we use more general structures following the Fox approach?
I wish you fruitful thoughts, and you can compare your ideas with that in my book, that has been 
in preparation for over 20 years, and whose few chapters are available in arXiv \cite{P-Book}. But here concisely:
\begin{enumerate}
\item [(1)] The oldest $n$-move conjecture is the Nakanishi 4-move conjecture: 
every knots can be reduced by 4-moves to the trivial knot. Formulated in 1979, it is still an open problem \cite{Kir}.\\
\item [(2)] One can look for an universal algebra (magma), $(X,*)$ where $*: X\times X \to X$ such that 
coloring of arcs of the diagram by elements of $X$ is consistent (Fig. 2.14) and is preserved by Reidemeister moves.
For example the third Reidemeister move leads to right self-distributivity $(a*b)*c=(a*c)*(b*c)$, Fig. 2.14. 
This leads to keis, racks, quandles, and shelfs as was explained in Scott Carter talk \cite{Ca}. The simplest example 
is ${\bf Z}_p$ with $a*b=2b-a$.
\end{enumerate}
\begin{center}
\psfig{figure=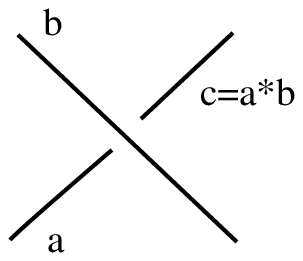,height=2.6cm},\ \ \ \ \ \ \ {} \psfig{figure=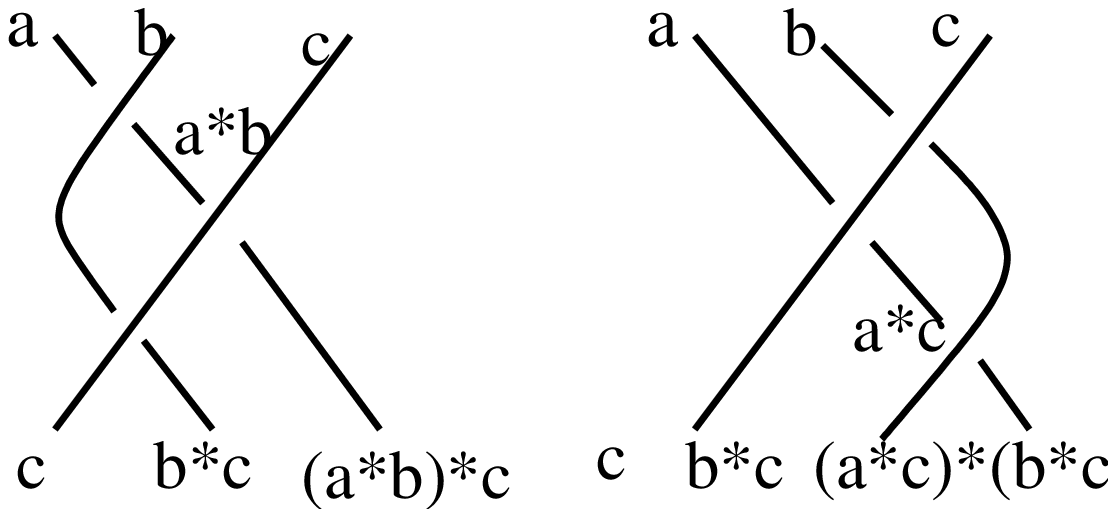,height=2.6cm}
\end{center}
\centerline{Fig. 2.14; coloring a crossing by elements of $X$ and the third Reidemeister move}

\section{Acknowledgement}
The author 
was partially supported by the NSA grant (\# H98230-08-1-0033), by the
Polish Scientific Grant: Nr. N-N201387034, by the GWU REF grant,  and by the
CCAS/UFF award.

\ \\
Department of Mathematics\\
George Washington University  \\
Washington, DC 20052 \\
USA\\
e-mail: przytyck@gwu.edu

\begin{thebibliography}{99 1000-2000}

\bibitem [A-B]{A-B}
J.~W.~Alexander, G.~B.~Briggs, On types of knotted curves,
{\it  Ann. of Math.}, 28(2), 1927/1928, 563-586.

\bibitem [A-V] {A-V}
J.~W.~Alexander, Letter to Oswald Veblen, 1919,
Papers of Oswald Veblen, 1881-1960 (bulk 1920-1960),
Archival Manuscript Material (Collection), Library of Congress.
(I would like to thank Jim Hoste for providing me with a copy
of the letter).

\bibitem[APR]{APR}
R.~P. Anstee, J.~H. Przytycki, D. Rolfsen. Knot polynomials and
generalized mutation.  {\it Topology and Applications}
32, 1989, 237-249;\\ 
 e-print:\ {\tt http://xxx.lanl.gov/abs/math.GT/0405382}

\bibitem[As] {As}
C.~W.~Ashley, {\it The Ashley book of knots}, Doubleday, 1944.

\bibitem [As] {As}
J.~Assion, Einige endliche Faktorgruppen der Zopfgruppen,
{\it Math. Z.}, 163, 1978, 291-302.

\bibitem [B-L]{B-L}
S.~Baldridge, A.~Lowrance,
Cube diagrams and a homology theory for knots;\\
e-print: \ {\tt arXiv:0811.0225v1 [math.GT]}

\bibitem[BLM]{BLM}
R.~D.~Brandt, W.~B.~R.Lickorish, K.~C.~Millett, A polynomial
invariant for unoriented knots and links, {\it Invent. Math. }
84, 1986, 563-573.

\bibitem [Brunn-1892]{Brunn-1892}
 H.~Brunn, Topologische Betrachtungen, {\it Zeitschrift f\"ur Mathematik
und Physik}, 37, 1892, 106-116.

\bibitem [B-Z]{B-Z}
G.~Burde, H.~Zieschang, {\it Knots}, De Gruyter, 1985 (second
edition 2003).

\bibitem[Ca]{Ca}
S.~Carter, A Survey of Quandle Ideas, This Book.

\bibitem[Che]{Che}
Q.Chen, The $3$-move conjecture for $5$-braids,
Knots in Hellas' 98;
The Proceedings of the International Conference
on Knot Theory and its Ramifications; Volume 1.
In the Series on Knots and Everything, Vol. 24, September 2000,
pp. 36-47.

\bibitem [Cr]{Cr}
R.H.Crowell, Knots and Wheels, National Council
of Teachers of Mathematics (N.C.T.M.) Yearbook, 1961.

\bibitem [C-F]{C-F}
R.H.Crowell, R.H.Fox, {\it An introduction to knot theory},
Ginn and Co., 1963.

\bibitem [D-P-1]{D-P-1}
M.~K.~D{\c a}bkowski, J.~H.~Przytycki, Burnside
obstructions to the Montesinos-Nakanishi 3-move conjecture, 
{\it Geometry and Topology (G\&T)}, 6, June, 2002, 335-360;\\
e-print:\ {\tt http://front.math.ucdavis.edu/math.GT/0205040}

\bibitem [De-0]{De-0}
M.~Dehn, {\it Jahresbericht der deutschen Mathematiker--Vereinigung},
Vol.16, 1907, p.573.

\bibitem[D-H] {D-H}
M.~Dehn, P.Heegaard, Analysis situs, {\it Encykl. Math. Wiss.}, vol. III
AB3 Leipzig, 1907, 153-220.

\bibitem[Del]{Del}
M.~Delbr\"uck, Knotting problems in biology, in {\it Mathematical Problems
in the Biological Sciences:} Proceedings of Symposia in Applied Mathematics, 14, 1962, 55-68.

\bibitem[Dia]{Dia}
Y.~Diao, Minimal knotted polygons on the cubic lattice,
{\it J. Knot Th. Ramifications}, {\bf 2}, 1993, 413-425.


\bibitem [DJP]{DJP}
J.~Dymara, T.~Januszkiewicz, J.~H.~Przytycki,
Symplectic structure on
Colorings,  Lagrangian tangles and Tits buildings, preprint, May 2001.

\bibitem [Ep]{Ep}
M.~Epple, Geometric aspects in the development of knot theory,
History of topology (ed. I.M.James) ,
301--357, North-Holland, Amsterdam, 1999.

\bibitem[F-R] {F-R}
R.~Fenn, C.~Rourke, Racks and links in codimension two,
{\it Journal of Knot Theory and its Ramifications}, 1(4) 1992, 343-406.


\bibitem [Fo-1]{Fo-1}
R.~H.~Fox, A quick trip through knot theory, {\it In: Top. 3-manifolds},
Proc. 1961 Top.Inst.Univ. Georgia (ed. M.K.Fort, jr), 120-167. Englewood
Cliffs. N.J.: Princeton-Hall, 1962.

\bibitem [Fo-2]{Fo-2}
R.~H.~Fox, Metacyclic invariants of knots and links, Canadian J.
Math., XXII(2) 1970, 193-201.

\bibitem[Gaus] {Gaus}
K.~F.~Gauss, Zur mathematischen Theorie der electrodynamischen  Wirkungen,
1833, {\it Werke}, K\"oniglichen Gesellschaft der Wissinchaften zu Gottingen,
5, 1877, 602-629.

\bibitem [Goe]{Goe}
L.~Goeritz, Knoten und quadratische Formen,
{\it Math. Z.}, 36,1933, 647-654.

\bibitem [Gor-1]{Gor-1}
C.~McA.~Gordon, Some aspects of classical knot theory,
{\it In: Knot theory}, L.N.M. 685, 1978, 1-60.

\bibitem[Gor-2]{Gor-2}
C.~McA.~Gordon, 3-Dimensional Topology up to 1960, 
{\it History of topology} (ed. I.M.James),
449--489, North-Holland, Amsterdam, 1999.

\bibitem[Heeg] {Heeg}
P.~Heegaard, Forstudier til en Topologisk Teori for de
algebraiske Fladers Sammenh{\ae}ng, K{\o}benhavn, 1898,
Filosofiske Doktorgrad;
French translation: Sur l'Analysis situs, {\em  Soc. Math. France Bull.},
44 (1916), 161-242. English translation, by A.Przybyszewska,
 of the topological part of Poul Heegaard Dissertation is in:
J.H.Przytycki, Knot theory from Vandermonde to Jones,
Preprint 43, Odense University 1993.

\bibitem[Ho] {Ho}
C.~F.~Ho, A new polynomial for knots and links; preliminary report,
{\it Abstracts AMS} 6(4), 1985, p. 300.

\bibitem [Joy]{Joy}
D.~Joyce, A classifying invariant of knots: the knot quandle,
{\it Jour. Pure Appl. Alg.}, 23, 1982, 37-65.

\bibitem[Kaw]{Kaw}
A.~Kawauchi, A survey of Knot Theory, Birkh\"ausen Verlag,
Basel-Boston-Berlin, 1996.

\bibitem[Kir]{Kir}
R.~Kirby, Problems in low-dimensional topology; Geometric Topology
(Proceedings of the Georgia International Topology Conference, 1993),
Studies in Advanced Mathematics, Volume 2 part 2., Ed. W.Kazez, AMS/IP,
1997, 35-473.

\bibitem[Klein]{Klein}
Felix Klein, \'Ueber den zusammenhang der fl\"achen. Mathematische Annalen, 9:476–482,\\
(some authors refer to this paper as containing a remark about no knotting in dimension four, on pages 476 or 478; 
however most likely the remark is not there.)

\bibitem [Li]{Li} W.B.R.Lickorish, A relationship between link polynomials,
{\it Math. Proc. Cambridge Phil. Soc.}, 100 (1986), 109-112.

\bibitem[Lis]{Lis}
J.~B.~Listing,  {\it  Vorstudien zur Topologie},  G\"ottinger Studien
(Abtheilung 1) 1, 1847,  811-875.

\bibitem[Mag]{Mag}
W.~Magnus, {\bf Max Dehn}, {\it Math. Intelligencer} 1 (1978),  132--143.
(Also in: Wilhelm Magnus Collected Paper, Edited by G.Baumslag and B.Chandler,
 Springer-Verlag, 1984).

\bibitem [Mat]{Mat}
S.~V.~Matveev,
Algorithmic Topology and Classification of 3-Manifolds {\it Algorithms and Computation in Mathematics}, Vol 9,
Second edition, Springer 2010 (First edition 2007).

\bibitem [Mo-1]{Mo-1}
J.~M.~Montesinos, Variedades de Seifert que son cubiertas ciclicas
ramificadas de dos hojas, {\it Bol. Soc. Mat. Mexicana (2)},
18, 1973, 1--32.

\bibitem [Mo-2]{Mo-2}
J.~M.~Montesinos, Lectures on $3$-fold simple coverings and $3$-manifolds,
Combinatorial methods in topology and algebraic geometry (Rochester,
N.Y., 1982),
Contemp. Math., 44, Amer. Math. Soc., Providence, RI, 1985, 157--177.

\bibitem [O'M]{O'M} 
O.~T.~O'Meara, Symplectic groups, Mathematics Surveys, Number 16, Published by AMS, 1978.

\bibitem [Pla]{Pla} A.~Plans, Contribution to the study of the homology groups of
the cyclic ramified coverings corresponding to a knot (Spanish),
{\it Revista Acad. Ci. Madrid} 47, 1953, 161--193. Math.Rev.
15,147a 56.0X (R.H.Fox).

\bibitem [Po-1]{Po-1}
 H.~Poincar\'e,  Analysis Situs (\&12), {\it  Journal d'Ecole Polytechnique
     Normale},   1, 1895,  1-121.

\bibitem [P-1]{P-1}
J.~H.~Przytycki, $t_k$-moves on links, In {\it Braids},
ed. J.S.Birman and A.Libgober, Contemporary Math. Vol. 78, 1988, 615-656.

\bibitem [P-2]{P-2}
J.~H.~Przytycki, Elementary conjectures in classical knot theory,
in {\em Quantum Topology}, Ed. L.J.Kauffman, R.A.Baadhio,
 Series on Knots and Everything - Vol.3, World Scientific, 1993, 292-320.

\bibitem[P-3]{P-3}
J.~H.~Przytycki, 3-coloring and other elementary
invariants of knots, Banach Center Publications, Vol. 42,
{\it Knot Theory}, 1998, 275-295; \\
e-print:\  {\tt  http://arxiv.org/abs/math.GT/0608172}


\bibitem[P-4]{P-4}
J.~H.~Przytycki, Classical roots of Knot Theory,
{\it Chaos, Solitons and Fractals}, Vol. 9 (No. 4-5), 1998, 531-545.

\bibitem[P-5]{P-5}
J.~H.~Przytycki,
From 3-moves to Lagrangian tangles and cubic skein modules,
{\it Advances in Topological Quantum Field Theory},
Proceedings of the NATO ARW on New Techniques in Topological Quantum
Field Theory, Kananaskis Village, Canada from 22 to 26 August 2001;
John M. Bryden (ed), October 2004, 71-125;\\
e-print: \ {\tt http://front.math.ucdavis.edu/math.GT/0405248}

\bibitem[P-6]{P-6}
J.~H.~Przytycki,
Kolorowanie splot\'ow (Coloring of knots), {\it Delta} 7, July 2003, 8-10.

\bibitem[P-Book]{P-Book}
J.~H.~Przytycki,
{\bf KNOTS:} From combinatorics of knot diagrams to the
combinatorial topology based on knots, Cambridge University Press,
accepted for publication, to appear 2013, pp. 600.\\
Chapter II, e-print:\ {\tt http://arxiv.org/abs/math/0703096}\\
Chapter V, e-print:\ {\tt http://arxiv.org/abs/math.GT/0601227} \\
Chapter IX, e-print:\ {\tt http://arxiv.org/abs/math.GT/0602264} \\
Chapter X, e-print:\ {\tt http://arxiv.org/abs/math.GT/0512630}\\

\bibitem [Rei-1]{Rei-1}
K.~Reidemeister,  Elementare Begrundung der Knotentheorie,
 {\it Abh. Math. Sem. Univ. Hamburg}, 5 (1927), 24-32.

\bibitem [Rol]{Rol}
D.~Rolfsen, {\it Knots and links}, Publish or Perish, 1976
(second edition, 1990).

\bibitem[Sq]{Sq}
C.~Squier, The Burau representation is unitary,
{\it Proc.  Amer. Math. Soc.}, 90(2), 1984, 199-202.

\bibitem [Uch] {Uch}
Y.~Uchida,
Proc. of Knots 96 ed. S. Suzuki 1997, 109-113
{\it Knots 96}, Proceedings of the Fifth International Research Institute
of MSJ, edited by Shin'ichi Suzuki, 1997 World Scientific Publishing Co.,
109-113.

\bibitem [Tait-1]{Tait-1}
P.~G.~Tait,  Preliminary note on a new method of investigating the properties 
of knots Proc. Royal Soc. Edinburgh, Vol. 9, 98 (1876-7), 403
(Tait writes here: "As we cannot have knots in two dimensions, and as Prof. Klein
has proved that they cannot exist in space of four dimensions;...") 


\bibitem [Tait-2]{Tait-2}
 P.~G.~Tait, On knots I,  II, III,  Scientific Papers,  Cambridge
 University Press, 1898-1900. Including: {\it Trans. R. Soc. Edin.},
28, 1877, 35-79.
Reprinted by Amphion Press, Washington D.C., 1993, (On Knots I, pages 273-317; remark 
about Klein is on the page 316).

\bibitem[Tur]{Tur}
V.~G.~Turaev, Quantum invariants of knots and 3-manifolds,
de Gruyter Studies in Mathematics, 18,
Walter de Gruyter \& Co., Berlin, 1994. x+588 pp.

\bibitem [Turing]{Turing}
A.~M.~Turing, Solvable and unsolvable problems, {Science News}, 31, 1954,
7-23.

\bibitem[Vi]{Vi}
O.~Ya.~Viro, Raskrasiennyje uzly, in Russian (Colored knots),
Kvant 3, 1981, 8-14 (Russian); English translation: Tied into
Knot Theory: unraveling the basics of
mathematical knots, Quantum 8(5), 1998, 16-20.

\bibitem[Wa-1]{Wa-1}
B.~Wajnryb, Markov classes in certain finite symplectic
representations of braid groups, in:  Braids (Santa Cruz, CA, 1986),
687--695, Contemp. Math., 78, Amer. Math. Soc.,
Providence, RI, 1988.

\bibitem[Wa-2]{Wa-2}
B.~Wajnryb, A braidlike presentation of ${\rm Sp}(n,p)$,
{\it Israel J. Math.} 76(3), 1991, 265--288.

\bibitem[Wir] {Wir}
W.~Wirtinger, \"Uber die Verzweigungen bei Funktionen von zwei
Ver\"anderlichen, {\it Jahresbericht d. Deutschen Mathematiker
Vereinigung}, 14 (1905), 517. (The title of the talk supposedly
given at September 26 1905 at the annual meeting of the German
Mathematical Society in Meran).

\bibitem [Wo-Kr]{Wo-Kr}
D.~Wolkstein, S.~N.~Kramer, Inanna: Queen of heaven and earth,
Harper Collins Publisher, 1983.

\end{thebibliography}
\end{document}